\documentclass[11pt]{article}
\usepackage{amsbsy,amsfonts,amsmath,amssymb,amsthm,euscript,enumerate}
\usepackage{graphicx,color,centernot,framed}
\definecolor{shadecolor}{gray}{0.85}
\usepackage{url}
\usepackage{hyperref}
\hypersetup{colorlinks,citecolor=blue,linkcolor=blue,urlcolor=black}
\usepackage{bbm}
\newcommand{\ind}[1]{\mathbbm{1}{\raisebox{-2pt}{$\scriptstyle \{#1\}$}}}
\newcommand{\indic}[1]{\mathbbm{1}{\raisebox{-2pt}{$\scriptstyle #1$}}}

\allowdisplaybreaks[1]

\theoremstyle{plain}
\numberwithin{equation}{section}
\theoremstyle{plain}
\newtheorem{thm}{Theorem}[section]

\newtheorem{lmm}[thm]{Lemma}
\newtheorem{prp}[thm]{Proposition}

\theoremstyle{definition}
\newtheorem{dfn}{Definition}
\newtheorem{rmk}{Remark}
\newtheorem*{AssD}{Assumption~D}
\newtheorem*{AssG}{Assumption~G}

\newcommand{\lbeq}[1]{\label{eq:#1}}

\newcommand{\nn}{\nonumber}

\newcommand{\Proof}[1]{\paragraph{\it #1}}
\newcommand{\QED}{\hspace*{\fill}\rule{7pt}{7pt}\bigskip}
\newcommand{\R}{{\mathbb R}}
\newcommand{\refe}[1]{(\ref{eq:#1})}

\newcommand{\sss}{\scriptscriptstyle}

\newcommand{\Z}{\mathbb{Z}}

\newcommand{\Zd}{\Z^d}
\newcommand{\Rd}{\R^d}
\newcommand{\norm}[1]{\lVert#1\rVert}

\makeatletter
\newcommand{\xRightarrow}[2][]{\ext@arrow 0359\Rightarrowfill@{#1}{#2}}
\makeatother

\definecolor{dgreen}{rgb}{0,0.6,0}

\title{Rate of convergence of the critical point of the memory-$\tau$ self-avoiding walk in dimensions $d>4$}
\author{
Noe Kawamoto 
\footnote{National Center for Theoretical Sciences, Taiwan. \url{https://orcid.org/0000-0003-2273-4497}}
}

\begin{document}
\maketitle
\begin{abstract}
We consider spread-out models of the self-avoiding walk and its finite-memory version, known as the memory-$\tau$ walk, which prohibits loops whose length is at most $\tau$, in dimensions $d>4$. The critical point is defined as the radius of convergence of the generating function for each model. It is known that the critical point of the memory-$\tau$ walk is non-decreasing in $\tau$ and converges to that of the self-avoiding walk as $\tau$ tends to infinity. In this paper, we study the rate at which the critical point of the memory-$\tau$ walk converges to that of the self-avoiding walk and show that the order is $\tau^{-(d-2)/2}$. The proof relies on the lace expansion, introduced by Brydges and Spencer. 
\end{abstract}


\section{Introduction}
The memory-$\tau$ walk is a model for walks that avoid loops of length $\tau$ or less. The self-avoiding walk is defined as the $\tau\rightarrow\infty$ limit of the memory-$\tau$ walk, representing an ensemble of walks with no self-intersections. These models are simply defined yet play a significant role in polymer chemistry and statistical physics. Moreover, both models are among the most important examples in the study of critical phenomena. For example,
the divergence of the generating function for the self-avoiding walk or the memory-$\tau$ walk starting at the origin at its radius of convergence is one of the best well-known examples of critical phenomena. This radius of convergence is referred to as the critical point, whose exact value is believed to depend on the details of the model. The lace expansion, introduced by Brydges and Spencer in \cite{bs85}, has been one of the most effective tools for analyzing  critical phenomena above the upper-critical dimension. Specifically, the lace expansion has played a crucial role in demonstrating mean-field behavior for various statistical mechanical models above their respective upper-critical dimensions. 

In this paper, we will study the critical points of the self-avoiding walk and the memory-$\tau$ walk in dimensions $d>4$ using the lace expansion. 

\subsection{Motivation and known results}
First, we define models. A path $\omega$ in $\Zd$ is defined as a sequence $\omega=(\omega(0), \omega(1),\cdots, \omega(|\omega|))$ of sites in $\Zd$, where $|\omega|$ denotes the number of steps of $\omega$. Let $D$ be a probability distribution on $\Zd$, which is invariant under all symmetries of $\Zd$. 
We define the weight function of $\omega$ by
\begin{align}\lbeq{Wpdef}
W_p(\omega)=\prod_{i=1}^{|\omega|}pD(\omega(i)-\omega(i-1)).
\end{align}  
By convention, when $|\omega|=0$, $W_p(\omega)$ is defined to be $1$.
Let $\mathcal W_n(x,y)$ be the set of the $n$-step walks with $\omega(0)=x$ and $\omega(n)=y$. In particular, $\mathcal W_0(x,y)$ consists of the $0$-step walk if $x=y$ and the empty set otherwise. To define walks that enforce self-avoidance constraint only within a time span $\tau$, we introduce for $1\le \tau \le \infty$
\begin{align}\lbeq{defust}
\mathcal U_{st}=-\delta_{\omega(s),\omega(t)},&&K_{\tau}[a,b]=\prod_{\substack{a \le s<t\le b \\ |s-t| \le \tau}}(1+\mathcal U_{st}).
\end{align}
We then define
\begin{align}
C^{\sss\tau}_{p,n}(x):=\sum_{\omega \in \mathcal W_n(o,x)}W_p(\omega)K_{\tau}[0,n].
\end{align}
This formulation simultaneously defines the self-avoiding walk ($\tau = \infty$) and the memory-$\tau$ walk ($\tau < \infty$). Throughout this paper, we omit the notation $\infty$ from any functions defined for the self-avoiding walk. The two-point function is defined by
\begin{align}
G_p^{\sss\tau}(x):=\sum_{n=0}^{\infty}C^{\tau}_{p,n}(x)=\sum_{\omega \in \mathcal W(o,x)}W_p(\omega)K_{\tau}[0,|\omega|],
\end{align}
where $\mathcal W(o,x)=\cup_{n=0}^{\infty}\mathcal W_n(o,x)$.
We define the susceptibility by the sum of the 2-point function  
\begin{align}\lbeq{chi-def}
\chi_p^{\sss\tau}:=\sum_{x\in\Zd}G_p^{\sss\tau}(x)=\sum_{n=0}^{\infty}p^nc_n^{\sss\tau},
\end{align}
where $\mathcal W_n=\cup_{x \in \Zd}\mathcal W_n(o,x)$, and 
\begin{align}\lbeq{cntau}
c_n^{\sss\tau}=\sum_{\omega \in \mathcal W_n}W_1(\omega)K_{\tau}[0,|\omega|].
\end{align}
The submultiplicative property $c_{n+m}^{\sss\tau}\le c_{n}^{\sss\tau}c_{m}^{\sss\tau}$, which follows from neglecting mutual avoidance, implies the existence of the connective constant 
\begin{align}
\mu_{\tau}=\lim_{n \rightarrow \infty}(c_{n}^{\sss\tau})^{1/n}=\inf_{n\ge1}(c_{n}^{\sss\tau})^{1/n}.
\end{align}
Therefore, the susceptibility has the radius of convergence given by
\begin{align}
p_c^{\sss\tau}:=\frac{1}{\mu_{\sss\tau}}
\end{align}
which we refer to as the critical point. The critical point is believed to be model-dependent, and its precise value has been studied for various statistical-mechanical models, particularly above their respective upper-critical dimension \cite{hs93, hs05, ks22, ms11, ms13}. 

For example, when $D$ follows a uniform distribution over $\{x\in\Zd:0<\max_j|x_j| \le L\}$ where $L$ is a sufficiently large positive constant, van der Hofstad and Sakai \cite{hs05} provide an explicit expression for $p_c$ of the self-avoiding walk in dimensions $d>4$: 
\begin{align}\lbeq{hs05pcinfty}
p_c=1+L^{-d}\sum_{n=2}^\infty U^{*n}(o)+O(L^{-d-1}),
\end{align}
where $U^{*n}$ denotes the $n$-fold convolution of $U$. The function $U$ represents the uniform distribution on the $d$-dimensional box $\{x\in\R^d:\max_j|x_j| \le1\}$, explicitly given by
\begin{align}\lbeq{defU}
U(x)=\frac{\ind{\norm{x}_{\infty}\le 1}}{2^d}.
\end{align}

Since $c_{n}^{\sss\tau}$ is decreasing in $\tau$, it follows that $\mu_{\sss\tau}$ is also decreasing in $\tau$. Thus, $\mu_{\sss\tau}$ always serves as an upper bound for $\mu$ for any $\tau$. Consequently, as $\tau \rightarrow \infty$, we have $p_c^{\sss\tau} \nearrow p_c$
as shown in \cite[Lemma~1.2.3]{ms93}. 

Various approaches are known for estimating $\mu_{\tau}$ in the nearest-neighbor model, where $D$ is given by  $D(x)=\frac1{2d}\ind{\norm{x}_1=1}$. For instance, in \cite[Appendix~A]{cfs59}, Fisher and Sykes discuss a method for determining the value of $\mu_{\sss\tau}$ by constructing a recurrence relation for $c_{n}^{\sss\tau}$. In their method, they classify a set of $n$-step memory-$\tau$ walks based on the number of steps required to close a loop of length $\tau$ or less. By analyzing the effect of adding an additional step to a path in each subset of $n$-step memory-$\tau$ walks, they derive a recurrence relation for the size of each subset. Then, $\mu_{\sss\tau}$ is given by the largest eigenvalue of the coefficient matrix of the recurrence relation. As the memory size increases, computing the eigenvalue becomes less practical. However, they mention that this method allows for the calculation of $\mu_{\sss\tau}$ up to $\tau=12$. In \cite{n98}, Noonan presents an alternative method for estimating $\mu_{\sss\tau}$. By partitioning a class of loops of length at most $\tau$ into subclasses, he derives an explicit generating function for the memory-$\tau$ walk and estimates $\mu_{\sss\tau}$ as the smallest positive root of the denominator of this generating function. In the paper, Noonan explicitly provides the generating function for the memory-$\tau$ for $\tau \le 8$. 

Kesten proved in \cite[Theorem~1]{k64} that there exists a constant $C_{d,\tau}$, depending on $d$ and $\tau$, such that for each even integer $\tau$ and $d\ge \max\{13\tau-13,5\}$,
\begin{align}\lbeq{kesten1}
\mu_{\sss\tau}-\mu\le C_{d,\tau}d^{-\frac{\tau+2}{2}}
\end{align}
where, with a constant $C_0$ independent of $d$ and $\tau$, 
\begin{align}
C_{d,\tau}=C_0\bigg[\frac{(\tau/2+2)^{\tau+4}}{(\tau/2+2)!}+\frac{1}{d}\frac{d!(d-4\tau+4)!}{[(d-2\tau+2)!]^2}\frac{(\tau/2+3)^{\tau+10}}{(\tau/2+3)!}\bigg].
\end{align}
Here, we note that \refe{kesten1} requires a condition on the relationship between $\tau$ and $d$, specifically that $d$ must be taken greater than $\max\{13\tau-13,5\}$. Thus, we can not take the limit $\tau \rightarrow \infty$ for fixed $d$ in \refe{kesten1}, whereas we can conclude from \refe{kesten1} that $\mu_{\sss\tau}-\mu=O(d^{-\frac{\tau+2}{2}})$ as $d\rightarrow \infty$ for fixed $\tau$. His approach to establishing \refe{kesten1} involves constructing a self-avoiding walk by removing loops of length at most $\tau$ from a memory-$\tau$ walk. By counting the number of memory-$\tau$ walks that turn into a single self-avoiding walk upon loop removal, he compares $c_n$ and $c_n^{\sss\tau}$ to examine the difference $c_n^{\sss\tau}-c_n$. The condition on the relationship between $\tau$ and $d$ stems from the requirement that two $\tau-1$-step memory-$\tau$ walks exist in completely disjoint subspaces with no shared dimensions. He combined \refe{kesten1} with the value of $\mu_4$, obtained using the approach outlined in \cite[Appendix~A]{cfs59} to derive
\begin{align}\lbeq{kesten2}
\mu=1-(2d)^{-1}-(2d)^{-2}+O(d^{-3}).
\end{align}
Furthermore, in \cite{hs93}, Hara and Slade showed that $\mu$ admits an asymptotic expansion in $(2d)^{-1}$ to all orders, with all coefficients being integers. Thanks to Kesten's result~\refe{kesten1}, it was sufficient for them to establish the existence of such a $(2d)^{-1}$- expansion of $\mu_{\tau}$. The proof is based on the lace expansion for the memory-$\tau$ walk. Moreover, they obtained a refined estimate compared to \refe{kesten2}: as $d \rightarrow \infty$
\begin{align}\lbeq{newmutau}
\mu&=1-(2d)^{-1}-(2d)^{-2}-3(2d)^{-3}-16(2d)^{-4}-102(2d)^{-5}+O(d^{-6}).
\end{align}
The proof of \refe{newmutau} relies on the lace expansion for the self-avoiding walk. They derive upper and lower bounds for the lace-expansion coefficient to compare their orders and identify the coefficients of terms with matching orders.

For the spread-out model, in which the interaction range $L$ is taken to be sufficiently large, Madras and Slade proved in \cite[Lemma~6.8.6]{ms93} that for $d>4$ and sufficiently large $L$, there exists a positive constant $K$ such that 
\begin{align}\lbeq{MadrasSlade}
p_c-p_c^{\sss\tau} \le K \tau^{-(1+\delta)}
\end{align}
where $\delta < (d-4)/2 \wedge 1$ and $K$ may depend on $d,~\delta,~L$ but not on $\tau$. We note that \refe{MadrasSlade} also holds for the nearest-neighbor model when $d$ is sufficiently large. The proof relies on the lace expansion for the memory-$\tau$ walk. Given two memories $\tau_1<\tau_2$, they obtain $\tau_2$-uniform upper bound of $p_c^{\sss\tau_2}-p_c^{\sss\tau_1}$ and then take the limit as $\tau_2 \rightarrow \infty$. They used the estimate \refe{MadrasSlade} to establish that there exists a constant $B$ such that $\sup_{x}C_{p_c,n}(x)\le Bn^{-d/2}$ \cite[Lemma~6.1.3]{ms93}, noting that this local central limit theorem-type result for the self-avoiding walk is derived through the analysis of the memory-$\tau$ walk. However, in \cite{hs02}, van der Hofstad and Slade successfully proved that there exists a positive constant $C$ such that $\sup_{x}C_{p_c,n}(x)\le CL^{-d}n^{-d/2}$ without relying on the memory-$\tau$ walk. 

Throughout this paper, we impose the following conditions on $D$, ensuring that $L^dD(Lx)$ serves as a discrete approximation of a function on $\mathbb R^d$. 
\begin{dfn}\label{dfn: D}
Let $h$ be a non-negative bounded function on $\Rd$ which is almost everywhere continuous, and which is invariant under reflections in coordinate hyperplanes and rotations by $\pi/2$. We suppose that there is an integrable function $H$ on $\Rd$ with $H(te)$ non-increasing in $t\ge 0$ for every unit vector $e\in \Rd$, such that $h(x)\le H(x)$. We assume that $\int_{\Rd}h(x)d^dx=1$ and $\int_{\Rd} |x|^{2+2\epsilon}h(x)d^dx <\infty$ for some $\epsilon >0$. Then, we define
\begin{align}\lbeq{defsigmah}
D(x)=\frac{h(x/L)}{\sum_{x \in \Zd}h(x/L)},&&\Sigma_h^2=\int |x|^2h(x)d^dx. 
\end{align}
\end{dfn}
\begin{rmk}
In \cite[Appendix~A]{hs02}, van der Hofstad and Slade confirm that the function $D$ defined above meets the following conditions, referred to as Assumption~D in \cite{hs02}. 
\begin{AssD}[\cite{hs02}]
For some $\epsilon>0$, $D$ has $2+2\epsilon$ moments, i.e.,
\begin{align}\lbeq{2plusepsmoment}
\sum_{x \in \Zd}|x|^{2+2\epsilon}D(x) < \infty,
\end{align}
and there exist constants $C_1, C_2, c_1, c_2$ and $\eta$ such that 
\begin{align}
\sigma^2:=\sum_{ x\in \Zd}{|x|}^2D(x)\le C_1L^2,~~~\sup_{x \in \Zd}D(x)\le C_2L^{-d},\lbeq{sigmabound}\\
c_1L^2|k|^2 \le 1-\hat D(k) \le c_2L^2|k|^2~~~(\norm{k}_{\infty} \le L^{-1}),\lbeq{LUL21-hatD}\\
1-\hat D(k) >\eta~~~(\norm{k}_{\infty}\ge L^{-1}),\lbeq{Leta1-hatD}\\
1-\hat D(k)<2-\eta~~~(k \in [-\pi,\pi]^d),\lbeq{U2-eta1-hatD}
\end{align}
where $\norm{k}_{\infty}=\max_i |k_i|$. 
\end{AssD}
\end{rmk}
A representative example of $D$ that satisfies all the assumptions is 
\begin{align}\lbeq{Ddef}
D(x)=\frac{\ind{0<\norm{x}_{\infty}\le L}}{(2L+1)^d-1},
\end{align}
for which $h(x)=2^{-d}\ind{0<\norm{x}_{\infty}\le 1}$.

\subsection{Main result}\label{Main result and outline of the proof}
The goal of this paper is to study the rate at which $p_c^{\sss \tau}$ converges to  $p_c$ as $\tau \rightarrow \infty$ when the range $L$ of $D$ is sufficiently large. 

Throughout this paper, $o_X(1)=A$ for $X=L$ or $X=\tau$ denotes that for any $\epsilon$, there exists $X_0$ such that for all $X\ge X_0$, we have $|A|\le \epsilon$. 

Our main result is stated in the following theorem.
\begin{thm}\label{thm:main}
Let $d>4$. For sufficiently large $L$ and $\tau$, we have
\begin{align}\lbeq{eqthmmain}
p_c-p_c^{\sss\tau}=\frac{2}{d-2}\left(\frac{d}{2\pi\Sigma_h^2}\right)^{\frac{d}{2}}L^{-d}\tau^{-\frac{d-2}{2}}\big[1+R(L,\tau)\big]
\end{align}
where  $\Sigma_h^2$ is defined by \refe{defsigmah} and the remainder term $R(L,\tau)$ satisfies
\begin{align}
|R(L,\tau)|\le o_L(1)+o_{\tau}(1).
\end{align}
Here, $R(L,\tau)$ is also dependent on $d$. 
\end{thm}

By Theorem~\ref{thm:main}, we conclude that the convergence rate of $p_c^{\sss\tau}$ to $p_c$ is of order $\tau^{-(d-2)/2}$. Moreover, we explicitly determine the proportionality constant for the leading term in terms of both $L$ and $\tau$. In these two aspects, \refe{eqthmmain} refines the previous result \refe{MadrasSlade}. The proof of Theorem~\ref{thm:main} relies on the lace expansions for both the self-avoiding walk and the memory-$\tau$ walk.

{\bf Notation.} 
\begin{itemize}
\item We use $C$ and $C'$ to denote positive constants that depend on $d$ but are independent of $\tau$ and $L$. Their values may change from line to line. 
\item The convolution of two absolutely summable functions $f,g:\Zd \rightarrow \mathbb C$ is defined for $x \in \Zd$ as
\begin{align}
(f*g)(x)=\sum_{y \in \Zd}f(y)g(x-y).
\end{align}
\item We use the Fourier transform and the inverse, defined for an absolutely summable functions $f:\Zd \rightarrow \mathbb C$ by
\begin{align}\lbeq{fourierdef}
\hat f(k)=\sum_{x \in \Zd}f(x)e^{i k\cdot x},~~~f(x)=\int_{[-\pi,\pi]^d}\hat f(k)e^{-ik\cdot x}\frac{d^dk}{(2\pi)^d},
\end{align}
where $k \in [-\pi,\pi]^d$ and $k\cdot x=\sum_{j=1}^{d}k_jx_j$. 

If the sum defining $\hat f$ is not well defined, we define $\hat f$ through the second identity of \refe{fourierdef}. The details on interpreting $\hat G_{p_c}$ are given in Remark~\ref{remark}. 

\item We denote the $L_{p}$-norm for $1\le p\le \infty$ by $\norm{\cdot}_{p}$ in both $x$-space and $k$-space. Specifically, for $f\in L_p(\Zd)$ and $\hat g \in L_p([-\pi,\pi]^d)$, we define
\begin{align}
&\norm{f}_{\infty}=\sup_{x \in \Zd}f(x),&&\norm{\hat g}_{\infty}=\sup_{k \in [-\pi,\pi]^d}\hat g(k),\nn\\
&\norm{f}_p=\bigg(\sum_{x\in \Zd}|f(x)|^p\bigg)^{1/p},&&\norm{\hat g}_{p}=\bigg(\int_{[-\pi,\pi]^d}|\hat g(k)|^p\frac{d^dk}{(2\pi)^d}\bigg)^{1/p}.
\end{align}
\end{itemize}

\subsection{The lace expansion}
\begin{figure}[t]
\begin{center}
\includegraphics[scale=0.6]{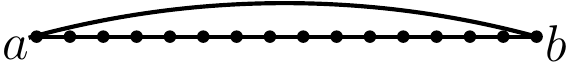}\\
\vskip2mm
\includegraphics[scale=0.6]{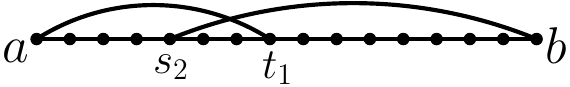}\\
\vskip2mm
\includegraphics[scale=0.6]{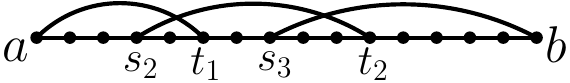}
\end{center}
\caption{Examples of $L \in \mathcal L_{\tau}^{\sss(N)}[a,b]$ for $N=1,2,3$. Each arc represents an edge, with edge lengths not exceeding $\tau$. When $b-a$ exceeds $\tau$, $\mathcal L_{\tau}^{\sss(1)}[a,b]=\varnothing$.}
\label{fig:Lace}
\end{figure}

Given an interval $[a,b]$ with nonnegative integers $a,b$, an open interval $(s,t)$  for $a\le s<t\le b$ is denoted by $st$ and referred to as an edge. We define a lace $L=\{s_1t_1,s_2t_2,\cdots,s_Nt_N\}$ on $[a,b]$ for $N\ge1$ as a set of $N$ edges satisfying the conditions: $s_1=a < s_2$, $s_{i+1}<t_i\le s_{i+2}$ for $1\le i\le N-2$ and $s_N<t_{N-1}<t_N=b$ (see Figure~\ref{fig:Lace}). 
We use  $\mathcal L_{\tau}^{\sss(N)}[a,b]$ to denote the set of laces on $[a,b]$,  each consisting of $N$ edges of length at most $\tau$. Recalling the definition of $\mathcal U_{st}$ given in \refe{defust}, we introduce
\begin{align}\lbeq{mathcalJN}
\mathcal J_{\tau}^{\sss(N)}[a,b]=\sum_{L\in \mathcal L_{\tau}^{\sss(N)}[a,b]}\prod_{st \in L}(-\mathcal U_{st})\prod_{s't' \in \mathcal C_{\tau}(L)}(1+\mathcal U_{s't'}),
\end{align}
where $\mathcal C_{\tau}(L)$ denotes the set of edges compatible with $L$ and of length at most $\tau$\footnote{Here, we define a compatible edge for a fixed $L$. A set of edges in $[a,b]$ is called a connected graph if $\bigcup_{st\in \Gamma}(s,t)=(a,b)$. Given a connected graph $\Gamma$, we define a lace $\mathsf L_{\Gamma}=\{s_1t_1, s_2t_2, \cdots\}$ associated with $\Gamma$ by determining $t_1,s_1,t_2,s_2,\cdots$ in the following way:  
\begin{enumerate}
\item $t_1=\max\{t: at \in \Gamma\}$,~~~$s_1=a$,
\item $t_{i+1}=\max\{t: \exists s<t_i~\text{such that}~st\in \Gamma\}$,~~~$s_{i+1}=\min\{s: st_{i+1} \in \Gamma\}$.
\end{enumerate} 
We iterate this procedure until $t_{i+1}=b$. Given a lace $L$, an edge $st \notin L$ is called a compatible edge with $L$  if it satisfies $\mathsf L_{L\cup \{st\}}=L$. For example, applying the procedure to the connected graph below yields the lace shown in Figure~\ref{fig:Lace} for $N=3$. 
\begin{center}
\includegraphics[scale=0.6]{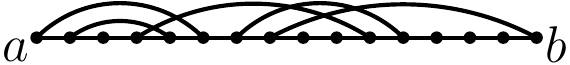}
\end{center}
}. Furthermore, we define 
\begin{align}\lbeq{pipnalpha}
\pi_{p, n}^{\sss\tau}(x)=\sum_{N=1}^{\infty}{(-1)}^N \pi_{p, n}^{\sss\tau, (N)}(x), && \pi_{p, n}^{\sss\tau, (N)}(x)=\sum_{\omega \in \mathcal W_n(o,x)}W_{p}(\omega)\mathcal J_{\tau}^{\sss(N)}[0,n].
\end{align}
The lace expansion gives the recursion equation for $n\ge0$ as
\begin{align}\lbeq{recursionGp}
C^{\sss\tau}_{p,n+1}(x)=p(D*C^{\sss\tau}_{p,n})(x)+\sum_{m=2}^{n+1}(\pi_{p, m}^{\sss\tau}*C^{\sss\tau}_{p,n+1-m})(x).
\end{align}
For the derivation of the lace expansion, see, for example, \cite{bs85, ms93, g004}. 

By applying the Fourier transform to both sides of \refe{recursionGp} for $k=0$ and summing over $n\ge 0$, and then solving for $\chi_p^{\sss\tau}$, we formally obtain
\begin{align}\lbeq{pc1-pi}
\chi_p^{\sss\tau}=\frac{1}{1-p-\hat \Pi_{p}^{\sss\tau}(0)},
\end{align}
where 
\begin{align}
\Pi_{p}^{\sss\tau,(N)}(x)=\sum_{m=2}^{\infty}\pi_{p,m}^{\sss\tau,(N)}(x),&&\Pi_{p}^{\sss\tau}(x)=\sum_{N=1}^{\infty}{(-1)}^N \Pi_{p}^{\sss\tau,(N)}(x).
\end{align}
Throughout this paper, we simplify notaion by omitting $(0)$ in the Fourier transform of functions evaluated at $k=0$. For example, we write $\hat \Pi_{p}^{\sss\tau}(0)$ as $\hat \Pi_{p}^{\sss\tau}$. We sometimes refer to $\hat \Pi_{p}^{\sss\tau}$ as the lace expansion coefficient.

\subsection{Outline of the proof}
Here we outline the proof of Theorem~\ref{thm:main}. Throughout this paper, we assume that $L$ is sufficiently large and frequently use the notation
\begin{align}\lbeq{defbeta}
\beta=L^{-d}.
\end{align}
Since $\chi_p^{\sss\tau}$ diverges as $p\rightarrow p_c^{\sss\tau}$, we obtain the identity 
\begin{align}\lbeq{mainidentity}
p_c^{\sss\tau}=1-\hat \Pi_{p_c^{\sss\tau}}^{\sss\tau}.
\end{align}

To establish Theorem~\ref{thm:main}, we express $p_c-p_c^{\sss\tau}$ in the following form. By substituting the expression from \refe{mainidentity}, we obtain
\begin{align}\lbeq{pc-pctau1}
p_c-p_c^{\tau}&=(1-\hat \Pi_{p_c})-(1-\hat \Pi_{p_c^{\sss \tau}}^{\sss\tau})\nn\\
&=\hat \Pi_{p_c}^{\sss(1)}-\hat \Pi_{p_c^{\sss\tau}}^{\sss\tau,(1)}+\sum_{N=2}^{\infty}(-1)^N\left[\hat \Pi_{p_c^{\sss\tau}}^{\sss\tau,(N)}-\hat \Pi_{p_c}^{\sss(N)}\right].
\end{align}
Using the expression from \refe{pc-pctau1}, it suffices to estimate the difference between the lace expansion coefficients for both the self-avoiding walk and the memory-$\tau$ walk at their respective critical points. To simplify the estimation, we further rewrite the difference of the lace expansion coefficients separately for the cases $N=1$ and $N\ge 2$. 

We first consider the case $N\ge 2$. To align the parameters with $p_c^{\sss\tau}$, we add and subtract $\hat \Pi_{p_c^{\sss\tau}}^{\sss(N)}$, resulting in 
\begin{align}\lbeq{pi2t-pi2}
\hat \Pi_{p_c^{\sss\tau}}^{\sss\tau,(N)}-\hat \Pi_{p_c}^{\sss(N)}&=\hat \Pi_{p_c^{\sss\tau}}^{\sss\tau,(N)}-\hat \Pi_{p_c^{\sss\tau}}^{\sss(N)}-(p_c-p_c^{\sss\tau})\frac{\hat \Pi_{p_c}^{\sss(N)}-\hat \Pi_{p_c^{\sss\tau}}^{\sss(N)}}{p_c-p_c^{\sss\tau}}\nn\\
&=\hat \Pi_{p_c^{\sss\tau}}^{\sss\tau,(N)}-\hat \Pi_{p_c^{\sss\tau}}^{\sss(N)}-(p_c-p_c^{\sss\tau})\partial_p \hat \Pi_{p_{\sss *}}^{\sss(N)}
\end{align}
for some $p_{\sss*}=p_{\sss *}^{\sss(N)} \in (p_c^{\sss\tau}, p_c)$. In the second equality, we use the differentiability of $\hat \Pi_{p}^{\sss(N)}$. Van der Hofstad and Slade proved in \cite{hGs03}\cite{hs02} that for $d>4$ and for sufficiently large $L$, there exists a positive constant $C$ such that
\begin{align}\lbeq{piNinfty1}
\sum_{x \in \Zd}|x|^q\pi_{p, n}^{\sss(N)}(x)\le (C\beta)^N\sigma^qN^qn^{-\frac{d-q}{2}}~~~(p\le p_c,~N\ge1,~n\ge2;~q=0,2,4).
\end{align}
By the definition of the weight function in \refe{Wpdef}, we have $\pi_{p, n}^{\sss(N)}(x)=p^n\pi_{1, n}^{\sss(N)}(x)$, which implies that $\partial_p\hat \pi_{p, n}^{\sss(N)}(0)=p^{-1}n\hat \pi_{p, n}^{\sss(N)}(0)$. Here, we use that $1\le p_c$\footnote{We note that $\sum_{\omega \in \mathcal W_n(o,x)}W_1(\omega)$ equals the $n$-step transition probability from $o$ to $x$ for the simple random walk with $1$-step transition probability given by $D$. By \refe{cntau}, $c_n^{\sss\tau}\le \sum_{\omega \in \mathcal W_n}W_1(\omega)=1$. Hence, the radius of convergence of the susceptibility should exceed $1$.}, where $1$ is the critical point for the simple random walk. Thus, using \refe{piNinfty1}, we obtain $\partial_p\hat \Pi_p^{\sss(N)}\le (C'\beta)^N$ for $p\le p_c$ and $N\ge 1$. Consequently, by the mean-value theorem, there exists some $p_{\sss *}=p_{\sss *}^{\sss(N)} \in (p_c^{\sss\tau}, p_c)$ such that $\frac{\hat \Pi_{p_c}^{\sss(N)}-\hat \Pi_{p_c^{\sss\tau}}^{\sss(N)}}{p_c-p_c^{\sss\tau}}=\partial_p \hat \Pi_{p_*}^{\sss(N)}$, where we note that $p_{\sss *}$ varies with $N$.

Next we consider the case $N=1$. Recall the definition given in \refe{pipnalpha}. Since the length of any edge in $[0,n]$ for $n\le \tau$ is less than $\tau$, it follows $\mathcal J^{\sss(N)}[0,n]=\mathcal J_{\sss\tau}^{\sss(N)}[0,n]$ for $N\ge1$ and $n\le \tau$. Consequently,  for $N\ge1$ and $n\le\tau$, we have
\begin{align}\lbeq{pinpitau}
\pi_{p, n}^{\sss(N)}(x)=\pi_{p, n}^{\sss\tau, (N)}(x).
\end{align}
Moreover, since $\mathcal L_{\tau}^{\sss(1)}[0,n]=\varnothing$ when $n\ge \tau+1$, (see Figure~\ref{fig:Lace}), it follows $\pi_{p, n}^{\sss\tau, (1)}(x)=0$ for $n\ge \tau+1$ and $x \in \Zd$, which leads to 
\begin{align}
\hat \Pi_{p}^{\sss\tau, (1)}=\sum_{n=2}^{\tau}\hat \pi_{p, n}^{\sss\tau, (1)}.
\end{align}
Using \refe{piNinfty1} again, we obtain for some $p_{\sss*}^{\sss(1)} \in (p_c^{\sss\tau}, p_c)$
\begin{align}\lbeq{decomppi1}
\hat \Pi_{p_c}^{\sss(1)}-\hat \Pi_{p_c^{\sss\tau}}^{\sss\tau,(1)}&=\sum_{n=2}^{\tau}\hat \pi_{p_c, n}^{\sss(1)}-\sum_{n=2}^{\tau}\hat \pi_{p_c^{\sss\tau}, n}^{\sss\tau, (1)}+\sum_{n=\tau+1}^{\infty}\hat \pi_{p_c, n}^{\sss(1)}\nn\\
&=(p_c-p_c^{\sss\tau})\sum_{n=2}^{\tau}\partial_p\hat \pi_{p_{\sss *}, n}^{\sss(1)}+\sum_{n=\tau+1}^{\infty}\hat \pi_{p_c, n}^{\sss(1)}.
\end{align}

Combining \refe{pc-pctau1}--\refe{pi2t-pi2} and \refe{decomppi1} yields 
\begin{align}\lbeq{pctau-pcinf}
p_c-p_c^{\sss\tau}&=K_{\sss\tau,L}\left\{\sum_{n=\tau+1}^{\infty}\hat \pi_{p_c, n}^{\sss(1)}+\sum_{N=2}^{\infty}(-1)^N\left[\hat \Pi_{p_c^{\sss\tau}}^{\sss\tau,(N)}-\hat \Pi_{p_c^{\sss\tau}}^{\sss(N)}\right]\right\}
\end{align}
where $K_{\sss\tau,L}$ is used to denote a $\tau$ and $L$-dependent constant, which is 
\begin{align}
K_{\sss\tau,L}=\left[1-\sum_{n=2}^{\tau}\partial_p\hat \pi_{p_{\sss*}, n}^{\sss(1)}+\sum_{N=2}^{\infty}(-1)^N\partial_p\hat \Pi_{p_{\sss*}}^{\sss(N)}\right]^{-1}.
\end{align}
Here, we observe that as indicated by \refe{piNinfty1}, $\partial_p\hat \pi_{p_{\sss*}, n}^{\sss(1)}\le C\beta n^{-d/2}$ and $\partial_p\hat \Pi_{p_{\sss*}}^{\sss(N)}\le (C'\beta)^N$, so there exist positive constants $c$ and $c'$ such that $1-c\beta^2\le K_{\sss\tau,L}\le 1+c'\beta$. 

Throughout this paper, we prove the following two propositions, which together establish Theorem~\ref{thm:main} with \refe{pctau-pcinf}. 
\begin{prp}\label{lmm:error0}
Let $d>4$. For sufficiently large $L$ and $\tau$,
\begin{align}
\sum_{n=\tau+1}^{\infty}\hat \pi_{p_c, n}^{\sss(1)}=\frac{2}{d-2}\left(\frac{d}{2\pi\Sigma_h^2}\right)^{\frac{d}{2}}\beta\tau^{-\frac{d-2}{2}}\big[1+R(L,\tau)\big]
\end{align}
with $|R(L,\tau)|\le o_L(1)+o_{\tau}(1)$. The error $R(L,\tau)$ is also dependent on $d$. 
\end{prp}
\begin{prp}\label{lmm:error1}
Let $d>4$. For sufficiently large $L$ and $\tau$, there exist positive constants $C$ and $C'$ such that
\begin{align}\lbeq{eqoflmmerror1}
C\beta^2\tau^{-\frac{d-2}{2}}\le \sum_{N=2}^{\infty}{(-1)}^N\left[\hat \Pi_{p_c^{\sss\tau}}^{\sss(N)}
-\hat \Pi_{p_c^{\tau}}^{\sss\tau,(N)}\right]\le C'\beta^2\tau^{-\frac{d-2}{2}}.
\end{align}
In particular, we obtain
\begin{align}
 \sum_{N=2}^{\infty}{(-1)}^N\left[\hat \Pi_{p_c^{\tau}}^{\sss(N)}
-\hat \Pi_{p_c^{\sss\tau}}^{\sss\tau,(N)}\right]=C_{\sss \tau, L}\beta^2\tau^{-\frac{d-2}{2}}.
\end{align}
where $C_{\sss\tau,L}$ is a constant that may depend on $\tau$ and $L$, and remains bounded as $\tau \rightarrow \infty$ or $L\rightarrow \infty$. 
\end{prp}
From the above two propositions, we see that both terms in \refe{pctau-pcinf} serve as the leading term with respect to $\tau$ whereas, in terms of $L$, the second term acts as error term. 

We prove Proposition~\ref{lmm:error0} in Section~\ref{Mainsection}, and Proposition~\ref{lmm:error1} in Section~\ref{Pi2proof}, respectively.

\subsection{Organization of the paper}
The remainder of the paper is organized as follows. In Section~\ref{results}, we introduce the results of van der Hofstad and Slade \cite{hs02}, which play a crucial role in our proof. Section~\ref{Mainsection} is devoted to proving  Proposition~\ref{lmm:error0}, which provides the leading term in terms of both $\tau$ and $L$. In Section~\ref{Pi2proof}, we establish Proposition~\ref{lmm:error1}. Finally, in Section~\ref{notes}, we discuss possible extensions of our result to other models. 

\section{The results by the inductive approach}\label{results}
We will heavily use the result achieved by van der Hofstad and Slade in \cite{hs02}. They established conditions under which the solution to the recursion relation derived from the lace expansion,
\begin{align}\lbeq{fgerecursion}
f_{n+1}(k;p)=\sum_{m=1}^{n+1}g_m(k;p)f_{n+1-m}(k;p)+e_{n+1}(k;p)~~~(n\ge0)
\end{align}
exhibits Gaussian behavior both in Fourier space and in the context of a local central limit theorem. Here, the functions $g_m$ and $e_m$ are given and the solution of \refe{fgerecursion} is $f_n(k;p)$, where $k \in [-\pi,\pi]^d$ and $p \ge 0$ are parameters. As for the self-avoiding walk, the Fourier transform of \refe{recursionGp}: 
\begin{align}\lbeq{FTrecursionGp}
\hat C_{p,n+1}(k)=p\hat D(k)\hat C_{p,n}(k)+\sum_{m=2}^{n+1}\hat\pi_{p, m}(k)\hat C_{p,n+1-m}(k)
\end{align}
is in the form of \refe{fgerecursion} with 
\begin{align}
&f_m(k;p)=\hat C_{p,m}(k)~~(m\ge0),&&g_1(k;p)=p\hat D(k),\nn\\
&g_m(k;p)=\hat\pi_{p, m}(k)~~(m\ge2),&&e_m(k;p)=0~~(m\ge 1).
\end{align} 
They introduce four conditions, Assumptions~S,D,G and E, that ensure the Gaussian behavior of the solution to \refe{fgerecursion}. Assumption~S concerns the symmetries of the functions in \refe{fgerecursion} with respect to the underlying lattice symmetries, as well as the uniform boundedness in $k \in [-\pi,\pi]^d$ of $f_n$ for each $n$. It is straightforward to verify that Assumption~S holds for the self-avoiding walk. Since $e_n\equiv 0$ for $n\ge1$ in the case of the self-avoiding walk, Assumption~E does not need to be checked. Assumption~D is given in Section~\ref{Main result and outline of the proof}. In \cite[Proposition~4.1]{hGs03}, van der Hofstad and Slade verified that Assumption G holds for the self-avoiding walk. Assumptions~G is as follows. 
\begin{AssG}[\cite{hs02}]
There is an $L_0$, an interval $I \subset [1-\delta,1+\delta]$ with $\delta \in (0,1)$ and a function $K_f\mapsto C_g(K_f)$, such that if the bounds 
\begin{align}\lbeq{127hs03}
\norm{{\hat D}^2f_m(\cdot;p)}_1\le K_f\beta m^{-\frac{d}{2}},~~~|f_m(0;p)|\le K_f,~~~|\nabla^2f_m(0;p)|\le K_f\sigma^2m
\end{align}
hold for some $K_f>1$, $L\ge L_0$, $p\in I$ and for all $1\le m\le n$, then 
for that $L$ and $p$, and for all $k \in [-\pi,\pi]^d$ and $2\le m\le n+1$, the following bounds hold:
\begin{align}
&\hskip15mm |g_m(k;p)|\le C_g(K_f)\beta m^{-\frac{d}{2}},~~~~|\nabla^2 g_m(0;p)|\le C_g(K_f)\sigma^2\beta m^{-\frac{d-2}{2}}, \lbeq{G1}\\
&\hskip30mm|\partial_p g_m(0;p)|=mp^{-1}g_m(0;p)\le C_g(K_f)\beta m^{-\frac{d-2}{2}},\lbeq{G2}\\
&\hskip20mm |g_m(k;p)-g_m(0;p)-[1-\hat D(k)]\sigma^{-2}\nabla^2g_m(0;p)|\nn\\
&\hskip30mm\le C_g(K_f)\beta [1-\hat D(k)]^{1+\epsilon'}m^{-\frac{d-2-2\epsilon'}{2}}\lbeq{G3}
\end{align}
where $\epsilon' \in [0,\epsilon]$ with $\epsilon$ given in \refe{2plusepsmoment}.
\end{AssG}
Under the four assumptions, \refe{127hs03} are shown to hold for all $m\ge1$ by induction. Therefore, all bounds \refe{G1}--\refe{G3} hold for all $m\ge2$ for models that satisfy the assumptions, including the self-avoiding walk. 

Before concluding this section, we summarize the results from \cite{hs02} that are needed for our analysis. 
\begin{thm}[$\cite{hs02}$]\label{thm:local}
Let $d>4$ and $L$ be sufficiently large. For $n \ge1$, the following hold:
\begin{itemize}
\item{Theorem 1.1(a)(d):} Fix $\gamma \in (0,1\wedge \frac{d-4}{4} \wedge \epsilon)$ and $\delta \in (0, (1\wedge \frac{d-4}{4} \wedge \epsilon )-\gamma)$, where $\epsilon$ is defined in \refe{2plusepsmoment}. Then, 
\begin{align}\lbeq{localCLT}
\hat C_{p_c,n}\Big(\frac{k}{\sqrt{v\sigma^2n}}\Big)=e^{-\frac{k^2}{2d}}\left[1+O(|k|^2n^{-\delta})+O(\beta n^{-\frac{d-4}{2}})\right]
\end{align}
where
\begin{align}\lbeq{vdef}
v=\frac{\sigma^2+p_c^{-1}\sum_{x\in \Zd}|x|^2\Pi_{p_c}(x)}{\sigma^2(1+\partial_p\hat \Pi_{p_c})}.
\end{align}
The error term is estimated uniformly in $\{k \in \mathbb R^d: 1-\hat D(k/\sqrt{v\sigma^2n})\le\gamma n^{-1}\log n\}$.
\item{(H4):} For $k \in \{k \in [-\pi,\pi]^d: 1-\hat D(k)> \gamma n^{-1}\log n\}$,
\begin{align}
\big|\hat C_{p_c,n}(k)\big|&\le C'_1(1-\hat D(k))^{-2-\rho}n^{-\frac{d}{2}},\lbeq{bCL1}\\
\big|\hat C_{p_c,n}(k)-\hat C_{p_c,n-1}(k)\big|&\le C'_2(1-\hat D(k))^{-1-\rho}n^{-\frac{d}{2}}.\lbeq{bCL2}
\end{align}
where $C'_1,C'_2$ are independent of $L$ and satisfy that $C'_1\gg C'_2$.
\end{itemize}
In the above, $\gamma, \delta, \rho>0$ are fixed such that $0<\frac{d-4}{2}-\rho<\gamma<\gamma +\delta<1\wedge \frac{d-4}{2}\wedge \epsilon$. 
\end{thm}
Furthermore, by \refe{G1} and \refe{G3} with $g_m(k;p)=\hat \pi_{p,m}(k)$, we obtain that for $p\le p_c$, there exist a positive constants $C'$ such that for all $k \in [-\pi,\pi]^d$,
\begin{align}
&\hskip15mm\big|\hat \Pi_{p}\big|\le C'\beta,~~~\bigg|\sum_{x\in \Zd}|x|^2\Pi_{p}(x)\bigg|\le C'\beta \sigma^2,\lbeq{Pidiffusion}\\ 
&\bigg|\hat \Pi_{p}-\hat \Pi_{p}(k)+(1-\hat D(k))\sigma^{-2}\sum_{x\in \Zd}|x|^2\Pi_{p}(x)\bigg|\le C'\beta (1-\hat D(k)).\lbeq{eee}
\end{align}
Here, we use $\nabla^2 g_m(0;p)=-\sum_{x}|x|^2\pi_{p,m}(x)$, which follows from the invariance of $\pi_{p,m}(x)$ under  translations and permutations, and \refe{piNinfty1}.

\section{Estimation of the proportionality constant}\label{Mainsection}
In this section, we prove Proposition~\ref{lmm:error0}. We begin by rewriting the expression that we aim to estimate. 
 
Since $\mathcal J^{\sss(1)}[0,n]=-\mathcal U_{0,n}\prod_{\substack{0\le s'<t'\le n\\ s't'\neq 0n}}(1+\mathcal U_{s't'})$, which equals $1$ if there are no self-intersections except for $\omega(0)=\omega(n)=o$, and $0$ otherwise, we can rewrite $\hat \pi_{p_c, n}^{\sss(1)}$ as 
\begin{align}
\hat \pi_{p_c, n}^{\sss(1)}=\sum_{x \in \Zd}\sum_{\omega \in \mathcal W_n(o,x)}W_{p_c}(\omega)\mathcal J^{\sss(1)}[0,n]=\sum_{\omega \in \mathcal W_n(o,o)}W_{p_c}(\omega)\prod_{\substack{0\le s'<t'\le n\\ s't'\neq 0n}}(1+\mathcal U_{s't'}). 
\end{align}
We claim that when $\omega \in \mathcal W_n(o,o)$,
\begin{align}\lbeq{1Ust}
\prod_{\substack{0\le s'<t'\le n\\ s't'\neq 0n}}(1+\mathcal U_{s't'})=\prod_{0<t'<n}(1+\mathcal U_{0t'})\prod_{1\le s'<t'\le n}(1+\mathcal U_{s't'})=K[1,n].
\end{align}
Since $\omega(0)=\omega(n)$, we have $\prod_{0<t'<n}(1+\mathcal U_{0t'})=\prod_{1\le s'<n}(1+\mathcal U_{s'n})$, which is included in the second product in the middle term of \refe{1Ust}. Thus, recalling the definition of $K[a,b]$ given in \refe{defust}, we obtain \refe{1Ust}.

By \refe{1Ust}, we derive
\begin{align}\lbeq{Piinf1>tau}
\sum_{n=\tau+1}^{\infty}\hat \pi_{p_c, n}^{\sss(1)}=\sum_{n=\tau+1}^{\infty}\sum_{\omega \in \mathcal W_n(o,o)}W_{p_c}(\omega)K[1,n]&=p_c\sum_{n=\tau+1}^{\infty}(D*C_{p_c,n-1})(o)\nn\\
&=p_c\int_{[-\pi,\pi]^d}\hat D(k)\sum_{n=\tau}^{\infty}\hat C_{p_c,n}(k)\frac{d^dk}{(2\pi)^d}. 
\end{align}
The following lemma provides an alternative expression for $\sum_{n=\tau}^{\infty}\hat C_{p_c,n}(k)$.
\begin{lmm}\label{lmm:rewriteCbiggertau}
Let $d>4$ and $\tau \ge 2$. Then, for $p \in [0,p_c]$
\begin{align}\lbeq{rewriteCbiggertau}
\sum_{n=\tau}^{\infty}\hat C_{p,n}(k)=\hat G_{p}(k)\left[\hat C_{p,\tau}(k)-\hat \pi_{p,\tau}(k)+\hat E_{p,\tau}(k)\right]
\end{align}
where
\begin{align}\lbeq{Ehatdef}
\hat E_{p,\tau}(k)=\sum_{m=2}^{\tau-1}\sum_{n=\tau+1-m}^{\tau-1}\hat \pi_{p,m}(k)\hat C_{p,n}(k)+\sum_{m=\tau}^{\infty}\sum_{n=0}^{\tau-1}\hat \pi_{p,m}(k)\hat C_{p,n}(k).
\end{align}
\end{lmm}
\Proof{Proof of Lemma~\ref{lmm:rewriteCbiggertau}.}In this proof, for simplicity, we denote $\sum_{n=\tau}^{\infty}\hat C_{p,n}(k)$ by $\hat C_{p,\ge\tau}(k)$.   

Summing both sides of \refe{FTrecursionGp} over $n\ge\tau$ and $\hat C_{p,0}=1$ gives
\begin{align}\lbeq{Cbiggertauk1}
\left[1-p\hat D(k)\right]\hat C_{p,\ge\tau}(k)=\hat C_{p,\tau}(k)+\sum_{n=\tau+1}^{\infty}\hat \pi_{p,n}(k)+\sum_{n=\tau}^{\infty}\sum_{m=2}^{n}\hat \pi_{p,m}(k)\hat C_{p,n+1-m}(k).
\end{align}
The last term on the right-hand side in \refe{Cbiggertauk1} can be rewritten as
\begin{align}\lbeq{Cbiggertauk2}
\sum_{n=\tau}^{\infty}\sum_{m=2}^{n}\hat \pi_{p,m}(k)\hat C_{p,n+1-m}(k)&=\sum_{m=2}^{\infty}\sum_{n=\tau\vee m}^{\infty}\hat \pi_{p,m}(k)\hat C_{p,n+1-m}(k)\nn\\
&=\sum_{m=2}^{\tau-1}\sum_{n=\tau}^{\infty}\hat \pi_{p,m}(k)\hat C_{p,n+1-m}(k)+\sum_{m=\tau}^{\infty}\sum_{n=m}^{\infty}\hat \pi_{p,m}(k)\hat C_{p,n+1-m}(k).
\end{align}
The first term of \refe{Cbiggertauk2} can be decomposed as 
\begin{align}
\sum_{m=2}^{\tau-1}\sum_{n=\tau}^{\infty}\hat \pi_{p,m}(k)\hat C_{p,n+1-m}(k)&=\sum_{m=2}^{\tau-1}\sum_{n=\tau+1-m}^{\infty}\hat \pi_{p,m}(k)\hat C_{p,n}(k)\nn\\
&=\sum_{m=2}^{\tau-1}\sum_{n=\tau+1-m}^{\tau-1}\hat \pi_{p,m}(k)\hat C_{p,n}(k)+\hat C_{p,\ge\tau}(k)\sum_{m=2}^{\tau-1}\hat \pi_{p,m}(k)
\end{align}
Similarly, for the second term of \refe{Cbiggertauk2}, we obtain
\begin{align}
\sum_{m=\tau}^{\infty}\sum_{n=m}^{\infty}\hat \pi_{p,m}(k)\hat C_{p,n+1-m}(k)&=\sum_{m=\tau}^{\infty}\sum_{n=1}^{\infty}\hat \pi_{p,m}(k)\hat C_{p,n}(k)\nn\\
&=\sum_{m=\tau}^{\infty}\sum_{n=0}^{\tau-1}\hat \pi_{p,m}(k)\hat C_{p,n}(k)+\hat C_{p,\ge\tau}(k)\sum_{m=\tau}^{\infty}\hat \pi_{p,m}(k)-\sum_{m=\tau}^{\infty}\hat \pi_{p,m}(k).
\end{align}
Therefore, the right-hand side on \refe{Cbiggertauk2} simplifies to
\begin{align}\lbeq{Cbiggertauk3}
\sum_{m=2}^{\infty}\hat \pi_{p,m}(k)\hat C_{p,\ge\tau}(k)+\hat E_{p, \tau}(k)-\sum_{m=\tau}^{\infty}\hat \pi_{p,m}(k).
\end{align}
Combining \refe{Cbiggertauk1} and \refe{Cbiggertauk3}, we obtain
\begin{align}
\left[1-p\hat D(k)-\sum_{m=2}^{\infty}\hat \pi_{p,m}(k)\right]\hat C_{p,\ge\tau}(k)=\hat C_{p,\tau}(k)-\hat \pi_{p,\tau}(k)+\hat E_{p,\tau}(k).
\end{align}
Since $\hat G_{p}(k)^{-1}=1-p\hat D(k)-\sum_{m=2}^{\infty}\hat \pi_{p,m}(k)$, we complete the proof. \QED

By \refe{Piinf1>tau} and Lemma~\ref{lmm:rewriteCbiggertau}, we can decompose $\sum_{n=\tau+1}^{\infty}\hat \pi_{p_c, n}^{\sss(1)}$ as
\begin{align}\lbeq{M1M2M3}
\sum_{n=\tau+1}^{\infty}\hat \pi_{p_c, n}^{\sss(1)}=p_c\left(M_1-M_2+M_3\right),
\end{align}
where
\begin{align}
M_1&=\int_{[-\pi,\pi]^d}\hat D(k)\hat C_{p_c,\tau}(k)\hat G_{p_c}(k)\frac{d^dk}{(2\pi)^d},\lbeq{defM1}\\
M_2&=\int_{[-\pi,\pi]^d}\hat D(k)\hat \pi_{p_c,\tau}(k)\hat G_{p_c}(k)\frac{d^dk}{(2\pi)^d},\lbeq{defM2}\\
M_3&=\int_{[-\pi,\pi]^d}\hat D(k)\hat E_{p_c,\tau}(k)\hat G_{p_c}(k)\frac{d^dk}{(2\pi)^d}.\lbeq{defM3}
\end{align}
We will show the following two lemmas in Section~\ref{D*} and in Section~\ref{derivativePiandintersectionPi1} respectively, by which we can conclude Proposition~\ref{lmm:error0}.

\begin{lmm}\label{lmm:M1Stau}
Let $d>4$. For sufficiently large $L$ and $\tau$, 
\begin{align}
M_1=\frac{2}{d-2}\left(\frac{d}{2\pi\Sigma_h^2}\right)^{\frac{d}{2}}\beta\tau^{-\frac{d-2}{2}}\big[1+R(L,\tau)\big]
\end{align}
with $|R(L,\tau)|\le o_L(1)+o_{\tau}(1)$. Here, $R(L,\tau)$ is also dependent on $d$. 
\end{lmm}

\begin{lmm}\label{lmm:M2M3}
Let $d>4$. For sufficiently large $L$ and $\tau$, 
\begin{align}\lbeq{M2M3}
|M_2|\le C\beta^2\tau^{-\frac{d}{2}},&&|M_3|\le C\beta^{\frac{3}{2}}\tau^{-\frac{d-2}{2}}.
\end{align}
\end{lmm}
\begin{rmk}
$M_2$ serves as an error term since it is small with respect to both $\beta$ and $\tau$ owing to $\hat \pi_{p_c,\tau}(k)$. In contrast, $\hat E_{p_c,\tau}(k)$ involves summations over the number of steps, $m$ and $n$, so it is not small in terms of $\tau$; however, due to $\hat \pi_{p_c,\cdot}(k)$ in $\hat E_{p_c,\tau}(k)$, it remains small with respect to $\beta$. 
\end{rmk}

\begin{rmk}\label{remark}
Here, we provide some comments on the interpretation of $\hat G_{p_c}$. Since $G_{p}(x)$ is not summable at $p=p_c$, we cannot define $\hat G_{p_c}(k)$ by the first identity in \refe{fourierdef}. Thus, we interpret $\hat G_{p_c}(k)$ using the second identity in \refe{fourierdef} as follows. 

Applying the Fourier transform to both terms of \refe{recursionGp}, summing over $n\ge1$ and solving the equation for $\hat G_p$, we obtain for $p<p_c$
\begin{align}\lbeq{1hatJpc}
\hat G_p(k)=\frac{1}{1-\hat J_p(k)}
\end{align}
where 
\begin{align}\lbeq{defJp}
J_{p}(x)=p\left(D(x)+p^{-1}\Pi_{p}(x)\right).
\end{align}
Using the identity given in \refe{pc1-pi} and the mean-value theorem with $p'\in (p,p_c)$, $1-\hat J_{p}(k)$ can be expressed as
\begin{align}\lbeq{1-hatJrewrite}
1-\hat J_{p}(k)&=p[1-\hat D(k)]+\hat \Pi_p-\hat \Pi_p(k)+(p_c-p)[1+\partial_p \hat \Pi_{p'}]\nn\\
&\ge [1-\hat D(k)]\bigg\{p+\frac{\hat \Pi_p-\hat \Pi_p(k)}{1-\hat D(k)}\bigg\},
\end{align}
where the inequality holds because the last term in the first line is positive for sufficiently large $L$, due to $1-\hat J_p(0)\ge 0$. 
Therefore, by \refe{Pidiffusion}--\refe{eee} and \refe{1hatJpc}--\refe{1-hatJrewrite}, we obtain that there exists a constant $C$ such that for $p< p_c$
\begin{align}\lbeq{infragauss}
\big|\hat G_{p}(k)\big| &\le \frac{C}{1-\hat D(k)}
\end{align}
for all $k\in [-\pi,\pi]^d$. 

Since $\hat G_p$ is integrable for $p<p_c$ and $d>2$, as follows from the bounds \refe{infragauss} and \refe{LUL21-hatD}--\refe{Leta1-hatD}, we use the left-continuity of $G_p$ in $p$ along with the dominated convergence theorem to obtain
\begin{align}
G_{p_c}(x)=\lim_{p \uparrow p_c}G_p(x)&=\int_{[-\pi,\pi]^d}\lim_{p \uparrow p_c}\frac{1}{1-\hat J_p(k)}e^{-ik\cdot x}\frac{d^dk}{(2\pi)^d}.\nn
\end{align}
Since $\hat \pi_{p,m}^{\sss(N)}(k)$ is left-continuous for $p\le p_c$, it follows from \refe{piNinfty1} that $\hat \Pi_p(k)$ is left-continuous for $p\le p_c$, which in turn implies the left-continuity of $\hat J_p$ for $p\le p_c$. Therefore, we obtain
\begin{align}\lbeq{kernelG}
G_{p_c}(x)
&=\int_{[-\pi,\pi]^d}\frac{1}{1-\hat J_{p_c}(k)}e^{-ik\cdot x}\frac{d^dk}{(2\pi)^d}\nn\\
&=\int_{[-\pi,\pi]^d}e^{-ik\cdot x}\int_{0}^{\infty} e^{-t\left[1-\hat J_{p_c}(k)\right]} dt\frac{d^dk}{(2\pi)^d}.
\end{align}
Throughout this paper, we define $\hat G_{p_c}(k)$ as the $t$-integral in the integrand with respect to $k$ in \refe{kernelG}. 
\end{rmk}

\subsection{Proof of Lemma~\ref{lmm:M1Stau}}\label{D*}
In this section, we prove Lemma~\ref{lmm:M1Stau}. Let
\begin{align}
S_{n}&=\left\{ k\in [-\pi,\pi]^d : 1-\hat D(k) \le \gamma n^{-1}\log n \right\},\\
L_{n}&=\left\{ k\in [-\pi,\pi]^d : 1-\hat D(k) > \gamma n^{-1}\log n \right\},
\end{align}
where $\gamma$ is defined by Theorem~\ref{thm:local}. By partitioning the integral domain, we can express $M_1$ as
\begin{align}\lbeq{ddd}
M_1=\int_{S_{\tau}}\hat D(k)\hat C_{p_c,\tau}(k)\hat G_{p_c}(k)\frac{d^dk}{(2\pi)^d}+\int_{L_{\tau}}\hat D(k)\hat C_{p_c,\tau}(k)\hat G_{p_c}(k)\frac{d^dk}{(2\pi)^d}.
\end{align}
We aim to show the following lemma which establishes Lemma~\ref{lmm:M1Stau}. 
\begin{lmm}\label{lmm:StauLtaulmm}
Let $d>4$. We fix $\gamma$ and $\delta$ defined in Theorem~\ref{thm:local}. Then, for sufficiently large $L$ and $\tau $, 
\begin{align}\lbeq{M1Stau}
\int_{S_{\tau}}\hat D(k)\hat C_{p_c,\tau}(k)\hat G_{p_c}(k)\frac{d^dk}{(2\pi)^d}=\frac{2}{d-2}\left(\frac{d}{2\pi\Sigma_h^2}\right)^{\frac{d}{2}}\beta\tau^{-\frac{d-2}{2}}[1+o_L(1)]+R_{\sss S}(\beta,\tau)
\end{align}
where $|R_{\sss S}(\beta,\tau)|\le \beta\tau^{-\frac{d-2}{2}}[o_L(1)+o_{\tau}(1)]+C\beta\tau^{-\frac{d-2}{2}-(K_2\wedge\delta)}$, with  $K_2\in (0, \frac{\gamma\Sigma_h^2}{4dc_2})$, and $c_2$ is defined in \refe{LUL21-hatD}. Here, $R_{\sss S}(L,\tau)$ is also dependent on $d$. 

Similarly, we have
\begin{align}\lbeq{M1Ltau0}
\Big|\int_{L_{\tau}}\hat D(k)\hat C_{p_c,\tau}(k)\hat G_{p_c}(k)\frac{d^dk}{(2\pi)^d}\Big|\le C\beta^{\frac{3}{2}}\tau^{-\frac{d-2}{2}}+C'\beta\tau^{-\frac{d-2}{2}}(\log \tau)^{-1}. 
\end{align}
\end{lmm}
We prove \refe{M1Stau} in Section~\ref{lmmM1stauproof} and \refe{M1Ltau0} in Section~\ref{lmmM1Ltau0proof} respectively. 

\subsubsection{Proof of \refe{M1Stau}.}\label{lmmM1stauproof}
In this section, we show \refe{M1Stau}. To prove \refe{M1Stau}, we use the following lemma.
\begin{lmm}\label{lmm:error of diffusion0}
Let $d>4$. Recall that $J_p$ is defined in \refe{defJp}. We define 
\begin{align}
\sigma_J^2:=\sum_{x\in\Zd}|x|^2J_{p_c}(x).
\end{align}
Then, we obtain that as $L\rightarrow \infty$, 
\begin{align}\lbeq{approximatesigmaJ}
\bigg|\frac{\sigma^2_J}{L^2}-\Sigma_h^2\bigg|\rightarrow 0
\end{align}
where $\Sigma_h^2$ is defined in \refe{defsigmah}
\end{lmm}
\Proof{Proof of Lemma~\ref{lmm:error of diffusion0}.} Recalling the definition of $J_p$ given in \refe{defJp}, we observe that
\begin{align}\lbeq{defsigmaJ}
\sigma_J^2=p_c\bigg(\sigma^2+p_c^{-1}\sum_{x\in \Zd}|x|^2\Pi_{p_c}(x)\bigg).
\end{align}
Thus, by \refe{sigmabound}, \refe{mainidentity}, \refe{Pidiffusion}, and \refe{defsigmaJ}, we obtain
\begin{align}\lbeq{sigmaJ}
|\sigma_J^2-\sigma^2|\le CL^2\beta.
\end{align}

By \refe{defsigmah}, we have 
\begin{align}\lbeq{sigma2L2conv}
\frac{\sigma^2}{L^2}=\frac{L^{-2}}{\sum_{x \in \Zd}h(x/L)}\sum_{x \in \Zd}|x|^2h(x/L)&=\frac{1}{\beta\sum_{y \in L^{-1}\Zd}h(y)}\beta\sum_{y\in L^{-1}\Zd}|y|^2h(y)
\end{align}
which converges to $\Sigma_h^2$ as $L \rightarrow \infty$. This implies that
\begin{align}
\bigg|\frac{\sigma_J^2}{L^2}-\Sigma_h^2\bigg|\le \bigg|\frac{\sigma^2}{L^2}-\Sigma_h^2\bigg|+C\beta,
\end{align}
where the right-hand side converges to $0$ as $L\rightarrow \infty$. 
\QED
\begin{rmk}
When considering the specific case that $h(x)=2^{-d}\ind{0<\norm{x}_{\infty}\le 1}$, for which $D$ is defined by \refe{Ddef}, we can observe the rate at which $\sigma/L^2$ converges to the variance of $U$ defined in \refe{defU}. In particular, we have
\begin{align}
\frac{\sigma^2}{L^{2}}=\int_{\Rd}{|x|}^2U(x)d^dx+O(L^{-1}).
\end{align}
\end{rmk}

Additionally, we will use 
\begin{align}\lbeq{D*2beta}
D^{*2}(o)=\beta \int_{\Rd}h(x)^2d^dx+o(\beta)
\end{align}
which follows from \refe{defsigmah}. 

From now on, for simplicity, we occasionally omit the differential element, i.e. $\frac{d^dk}{(2\pi)^d}$. 
\Proof{Proof of \refe{M1Stau}.}
By \refe{localCLT} and \refe{kernelG}, the left-hand side of \refe{M1Stau} can be decomposed as
\begin{align}\lbeq{M11M12}
\int_{S_{\tau}}\hat D(k)\hat C_{p_c,\tau}(k)\hat G_{p_c}(k)=\left[1+O\big(\beta \tau^{-\frac{(d-4)}{2}}\big)\right]M_{1.1}+M_{1.2}
\end{align}
where
\begin{align}
M_{1.1}&=\int_{0}^{\infty}dt\int_{S_{\tau}}\hat D(k)e^{-\tau\frac{v\sigma^2|k|^2}{2d}}e^{-t\left[1-\hat J_{p_c}(k)\right]},\lbeq{M11}\\
M_{1.2}&=\int_{S_{\tau}}\hat D(k)e^{-\tau\frac{v\sigma^2|k|^2}{2d}}\hat G_{p_c}(k)O(L^2|k|^2\tau^{1-\delta}).\lbeq{M12}
\end{align}

We begin by analyzing \refe{M11}. Our goal is to establishing that there exist constnats $C, K_1$ and $K_2$ such that 
\begin{align}\lbeq{M11final}
&\Big|M_{1.1}-\frac{2}{d-2}\left(\frac{d}{2\pi\Sigma_h^2}\right)^{\frac{d}{2}}\beta\tau^{-\frac{d-2}{2}}[1+o_L(1)]\Big|\nn\\
&\le C\beta\tau^{-\frac{d-2}{2}}\left[e^{-K_1\tau L^2}+\tau^{-K_2}+\tau^{-1}\log \tau+o_{L,\tau}(1)\right].
\end{align}
Since $1-\hat D(k)\le \gamma \tau^{-1}\log \tau$ for $k \in S_{\tau}$, we obtain
\begin{align}\lbeq{M11'}
\Big|M_{1.1}-\int_{0}^{\infty}dt\int_{S_{\tau}}e^{-\tau\frac{v\sigma^2|k|^2}{2d}}e^{-t\left[1-\hat J_{p_c}(k)\right]}\Big|&\le \int_{0}^{\infty}dt\int_{S_{\tau}}|1-\hat D(k)|e^{-\tau\frac{v\sigma^2|k|^2}{2d}}e^{-t\left[1-\hat J_{p_c}(k)\right]}\nn\\
&\le \gamma \tau^{-1}\log \tau\int_{0}^{\infty}dt\int_{S_{\tau}}e^{-\tau\frac{v\sigma^2|k|^2}{2d}}e^{-t\left[1-\hat J_{p_c}(k)\right]}.
\end{align}
Therefore, to estimate $M_{1.1}$, it suffices to analyze the integral $\int_{0}^{\infty}dt\int_{S_{\tau}}e^{-\tau\frac{v\sigma^2|k|^2}{2d}}e^{-t\left[1-\hat J_{p_c}(k)\right]}$. 
After we prove \refe{M1Stau}, we will show that
\begin{align}\lbeq{M11c}
\int_{0}^{\infty}dt\int_{S_{\tau}}e^{-\tau\frac{v\sigma^2|k|^2}{2d}}e^{-t\left[1-\hat J_{p_c}(k)\right]}-\int_{0}^{\infty}dt\int_{S_{\tau}}e^{-(\tau+t)\frac{\sigma_J^2|k|^2}{2d}}=\beta\tau^{-\frac{d-2}{2}}o_{L\cap \tau}(1).
\end{align}
Here, $o_{L\cap \tau}(1)=A $ means that for any $\epsilon$, there exist a pair $(L_o, \tau_o)$ such that for all $ (L,\tau)\ge(L_o,\tau_o) $, $|A|\le \epsilon$ holds. The notation $(a,b)\ge (a_0,b_0)$ signifies that $a\ge a_0$ and $b\ge b_0$. 

Therefore, it suffices to show that there exist constants $C$, $K_1\in (0,\frac{\pi^2\Sigma_h^2}{4d})$ and $K_2\in (0, \frac{\gamma\Sigma_h^2}{4dc_2})$, such that 
\begin{align}\lbeq{mainI123}
&\left|\int_{0}^{\infty}dt\int_{S_{\tau}}e^{-(\tau+t)\frac{\sigma_J^2|k|^2}{2d}}-\frac{2}{d-2}\left(\frac{d}{2\pi\sigma_J^2}\right)^{\frac{d}{2}}\tau^{-\frac{d-2}{2}}\right|\nn\\
&\le C\beta\tau^{-\frac{d-2}{2}}\left[e^{-K_1\tau L^2}+\tau^{-K_2}\right]. 
\end{align}
By the change of variable of $t$, the first term on the left-hand side of \refe{mainI123} can be decomposed as
\begin{align}\lbeq{III123}
\int_{0}^{\infty}dt\int_{S_{\tau}}e^{-(\tau+t)\frac{\sigma_J^2|k|^2}{2d}}=\int_{\tau}^{\infty}dt\left(I_{t,1}-I_{t,2}-I_{t,3}\right)
\end{align}
where
\begin{align}
I_{t,1}=\int_{\mathbb R^d} e^{-t\frac{\sigma_J^2|k|^2}{2d}},&&I_{t,2}=\int_{\mathbb R^d\backslash [-\pi,\pi]^d} e^{-t\frac{\sigma_J^2|k|^2}{2d}},&&I_{t,3}=\int_{L_{\tau}} e^{-t\frac{\sigma_J^2|k|^2}{2d}}.
\end{align}
Using \refe{approximatesigmaJ}, the first integral $I_{t,1}$ is estimated exactly as
\begin{align}\lbeq{It1}
I_{t,1}=\left(\frac{d}{2\pi\sigma_J^2t}\right)^{\frac{d}{2}}=\left(\frac{d}{2\pi\Sigma_h^2}\right)^{\frac{d}{2}}\left[1+o_L(1)\right]\beta t^{-\frac{d}{2}}.
\end{align}
For $I_{t,2}$ and $I_{t,3}$, we apply 
\begin{align}\lbeq{techn}
\int_{|k|\ge b}e^{-a{|k|}^2}\frac{d^dk}{{(2\pi)}^d}\le \exp(-ab^2/2){\left(2\pi a\right)}^{-d/2}~(\forall a,b>0).
\end{align}
Since $t\ge\tau$, we can bound each of them as follows. 
\begin{align}
I_{t,2}&\le \int_{|k|> \pi} e^{-t\frac{\sigma_J^2|k|^2}{2d}}\le \left(\frac{d}{\pi\sigma_J^2t}\right)^{\frac{d}{2}}e^{-\tau\frac{\pi^2\sigma^2_J}{4d}}\le C\beta e^{-K_1\tau L^2}t^{-\frac{d}{2}},\\
I_{t,3}&\le \int_{|k|^2 \ge \gamma c_2^{-1} L^{-2}\tau^{-1}\log \tau}e^{-t\frac{\sigma_J^2|k|^2}{2d}}\le \left(\frac{d}{\pi\sigma_J^2t}\right)^{\frac{d}{2}}\tau^{-\frac{\gamma \sigma_J^2}{4dc_2L^2}}\le C\beta\tau^{-K_2}t^{-\frac{d}{2}}.\lbeq{It3}
\end{align}
Here, for the first inequality in \refe{It3}, we use $L_{\tau}\subset \{k\in [-\pi,\pi]^d: c_2L^2|k|^2> \gamma\tau^{-1}\log \tau \}$, which follows from \refe{LUL21-hatD}. Combining \refe{III123} and \refe{It1}--\refe{It3}, we obtain \refe{mainI123}.

Finally we consider \refe{M12}. By using \refe{infragauss} and performing the change of variables from $k$ to $l:=L\sqrt{\tau}k$, we can bound $M_{1.2}$ as
\begin{align}\lbeq{abM12}
|M_{1.2}|&\le \int_{S_{\tau}}e^{-\tau\frac{v\sigma^2|k|^2}{2d}}|\hat G_{p_c}(k)|O(L^2|k|^2\tau^{1-\delta})\frac{d^dk}{(2\pi)^d}\nn\\
&\le C\int_{S_{\tau}}e^{-\tau\frac{v\sigma^2|k|^2}{2d}}\frac{1}{1-\hat D(k)}O(L^2|k|^2\tau^{1-\delta})\frac{d^dk}{(2\pi)^d}\nn\\
&=C\beta\tau^{-\frac{d}{2}-\delta}\int_{\tilde S_{\tau}}e^{-\frac{v\sigma^2|l|^2}{2dL^2}}\frac{1}{1-\hat D(l/L\sqrt{\tau})}O(|l|^2)\frac{d^dl}{(2\pi)^d}
\end{align}
where $\tilde S_{\tau}$ denotes the domain of $l$ corresponding to $S_{\tau}$. Since by \refe{LUL21-hatD}--\refe{Leta1-hatD}, we obtain 
\begin{align}
\frac{1}{1-\hat D(l/(L\sqrt{\tau}))}\le
\begin{cases}
\frac{1}{c_1}\tau |l|^{-2}&(\norm{l}_{\infty}\le \sqrt{\tau}),\\
\eta^{-1}&(\norm{l}_{\infty}\ge \sqrt{\tau}),
\end{cases}
\end{align}
the integral of \refe{abM12} can be bounded above by a constant multiple of $\tau$. Therefore, we obtain
\begin{align}\lbeq{approxM12}
|M_{1.2}|\le C\beta\tau^{-\frac{d-2+2\delta}{2}}.
\end{align}
By combining \refe{M11M12}, \refe{M11final} and \refe{approxM12}, we complete the proof of \refe{M1Stau}.
\QED

\Proof{Proof of \refe{M11c}.}
We show \refe{M11c}. By the change of variable of $t$, the left-hand side of \refe{M11c} becomes 
\begin{align}
\int_{\tau}^{\infty}dt \int_{S_{\tau}}\left\{e^{-\tau\frac{v\sigma^2|k|^2}{2d}}e^{-(t-\tau)\left[1-\hat J_{p_c}(k)\right]}-e^{-t\frac{\sigma_J^2|k|^2}{2d}}\right\}.
\end{align}
By the change of variable from $k$ to $l:=L\sqrt{t}k$, we obtain
\begin{align}
\int_{S_{\tau}}\left\{e^{-\tau\frac{v\sigma^2|k|^2}{2d}}e^{-(t-\tau)\left[1-\hat J_{p_c}(k)\right]}-e^{-t\frac{\sigma_J^2|k|^2}{2d}}\right\}=\beta t^{-\frac{d}{2}}I_{t,4}
\end{align}
where
\begin{align}\lbeq{defI4}
I_{t,4}=\int_{\tilde S_{\tau}}\Big\{e^{-\frac{\tau v\sigma^2|l|^2}{2dL^{2}t}}e^{-(t-\tau)\left[1-\hat J_{p_c}\big(\frac{l}{L\sqrt{t}}\big)\right]}-e^{-\frac{\sigma_J^2|l|^2}{2dL^{2}}}\Big\}\frac{d^dl}{(2\pi)^d}.
\end{align}
Here, $\tilde S_{\tau}$ is used to denote the domain of $l$ corresponding to $S_{\tau}$:
\begin{align}\lbeq{Stau}
\tilde S_{\tau}=\{l \in \big[-L\sqrt{t}\pi,L\sqrt{t}\pi\big]^d: 1-\hat D\big(l/L\sqrt{t}\big)\le \gamma \tau^{-1}\log \tau\}.
\end{align}
To prove \refe{M11c}, it suffices to show that
\begin{align}\lbeq{d}
\int_{\tau}^{\infty}\beta t^{-\frac{d}{2}}I_{t,4}dt=\beta \tau^{-\frac{d-2}{2}}o_{L\cap \tau}(1).
\end{align}

Let $[\cdots]$ denote the integrand of \refe{defI4}. Since by \refe{LUL21-hatD}--\refe{Leta1-hatD}, \refe{1hatJpc} and \refe{infragauss}, we have
\begin{align}\lbeq{dd}
\Big|1-\hat J_{p_c}\Big(\frac{l}{L\sqrt{t}}\Big)\Big|=\Big|\hat G_{p_c}\Big(\frac{l}{L\sqrt{t}}\Big)\Big|^{-1}\ge C\Big[1-\hat D\Big(\frac{l}{L\sqrt{t}}\Big)\Big]\ge 
\begin{cases}
C|l|^2t^{-1}&[\norm{l}_{\infty}\le\sqrt{t}],\\
C\eta&[\norm{l}_{\infty}> \sqrt{t}],
\end{cases}
\end{align}
applying \refe{techn} together with \refe{Stau} and \refe{dd} yields
\begin{align}\lbeq{m}
\int_{\tilde S_{\tau},~ \norm{l}_{\infty} >\sqrt{t}}[\cdots]
&\le\int_{ \tilde S_{\tau},~ \norm{l}_{\infty} >\sqrt{t}}\Big\{e^{-\frac{\tau v \sigma^2 |l|^2}{2dL^2t}}e^{-(t-\tau)\big[1-\hat D\big(\frac{l}{L\sqrt{t}}\big)\big]}+e^{-\frac{\sigma_J^2 |l|^2}{2dL^2}}\Big\}\nn\\
&=\int_{ \tilde S_{\tau},~ \norm{l}_{\infty} >\sqrt{t}}\Big\{e^{-\frac{\tau v \sigma^2 |l|^2}{2dL^2t}}e^{-t\big[1-\hat D\big(\frac{l}{L\sqrt{t}}\big)\big]}e^{\tau\big[1-\hat D\big(\frac{l}{L\sqrt{t}}\big)\big]}+e^{-\frac{\sigma_J^2 |l|^2}{2dL^2}}\Big\}\nn\\
&\overset{\refe{Stau},\refe{dd}}{\le} \int_{ \tilde S_{\tau},~ \norm{l}_{\infty} >\sqrt{t}}\Big\{e^{-\frac{\tau v \sigma^2 |l|^2}{2dL^2t}}e^{-C\eta t}e^{\gamma\log\tau}+e^{-\frac{\sigma_J^2 |l|^2}{2dL^2}}\Big\}\nn\\
&\le e^{-C\eta t}e^{\gamma\log\tau}\Big(\frac{dL^2t}{\pi \tau v\sigma^2}\Big)^{\frac{d}{2}}e^{-\tau\frac{v\sigma^2}{2dL^2}}+\Big(\frac{dL^2}{\pi \sigma_J^2}\Big)^{\frac{d}{2}}e^{-\frac{\sigma_J^2 t}{2dL^2}}\nn\\
&\le C'e^{-C\tau}.
\end{align} 
Here, for the last inequality, we use the bound $v\sigma^2\ge \sigma_J^2$ (see \refe{vsig-sigJ}--\refe{PiderivativePi} below) and Lemma~\ref{lmm:error of diffusion0}. Therefore, we obtain
\begin{align}
\int_{\tau}^{\infty}\beta t^{-\frac{d}{2}}I_{t,4}dt
&\le \beta \int_{\tau}^{\infty}t^{-\frac{d}{2}}\int_{\tilde S_{\tau},~\norm{l}_{\infty} \le\sqrt{t}}[\cdots]d^dl+C'e^{-C\tau}\beta \int_{\tau}^{\infty}t^{-\frac{d}{2}}\nn\\
&\le \beta \int_{\tau}^{\infty}t^{-\frac{d}{2}}\int_{\tilde S_{\tau},~\norm{l}_{\infty} \le\sqrt{t}}[\cdots]d^dl+C''\beta e^{-C\tau}\tau^{-\frac{d-2}{2}}.
\end{align}
Thus, to prove \refe{d}, we aim to establish 
\begin{align}\lbeq{gg}
\beta\int_{\tau}^{\infty}t^{-\frac{d}{2}}\int_{\tilde S_{\tau},~\norm{l}_{\infty} \le\sqrt{t}}[\cdots]d^dl=\beta \tau^{-\frac{d-2}{2}}o_{L\cap \tau}(1).
\end{align}
To apply the dominated convergence theorem, we first show that the inegral part with respect to $l$ on the left-hand side of \refe{gg} can be bounded uniformly with repect to both $L$ and $t$. By \refe{dd}, we obtain
\begin{align}
\int_{\tilde S_{\tau},~ \norm{l}_{\infty} \le \sqrt{t}}[\cdots]
&\le \int_{\tilde S_{\tau},~\norm{l}_{\infty} \le \sqrt{t}}\Big|e^{-t\big[1-\hat J_{p_c}\big(\frac{l}{L\sqrt{t}}\big)\big]}e^{-\tau\frac{(v\sigma^2-\sigma_J^2)|l|^2}{2dL^2t}}e^{\tau\Big[1-\hat J_{p_c}\big(\frac{l}{L\sqrt{t}}\big)-\frac{\sigma_J^2|l|^2}{2dL^2t}\Big]}-e^{-\frac{\sigma_J^2 |l|^2}{2dL^2}}\Big|\nn\\
&\overset{\refe{dd}}{\le} \int_{\tilde S_{\tau},~\norm{l}_{\infty} \le \sqrt{t}}\Big|e^{-t\big[1-\hat J_{p_c}\big(\frac{l}{L\sqrt{t}}\big)\big]}e^{-\tau\frac{(v\sigma^2-\sigma_J^2)|l|^2}{2dL^2t}}-e^{-\frac{\sigma_J^2 |l|^2}{2dL^2}}\Big|\lbeq{ee}\\
&\le C'\int_{\mathbb R^d}e^{-C|l|^2}<\infty.\lbeq{cc}
\end{align}
Here, for the third and forth lines, we use $1-\hat J_{p_c}\Big(\frac{l}{L\sqrt{t}}\Big)-\frac{\sigma_J^2|l|^2}{2dL^{2}t}\le 0$ and $v \sigma^2-\sigma_J^2\ge 0$, respectively. Indeed, by the Taylor expansion:
\begin{align}\lbeq{taylor}
e^{ix}=1+ix-\frac{1}{2}x^2-x^2\int_{0}^{1}(1-s)(e^{isx}-1)ds~~~(x \in \mathbb R),
\end{align}
we obtain
\begin{align}\lbeq{358}
1-\hat J_{p_c}\bigg(\frac{l}{L\sqrt{t}}\bigg)-\frac{\sigma_J^2|l|^2}{2dL^{2}t}=\sum_{x\in \Zd}J_{p_c}(x)\frac{(l\cdot x)^2}{L^{2}t}\int_{0}^{1}(1-s)(e^{\frac{isl\cdot x}{L\sqrt{t}}}-1)\le 0.
\end{align}
Furthermore, by \refe{mainidentity}, \refe{vdef} and \refe{defsigmaJ}, we obtain
\begin{align}\lbeq{vsig-sigJ}
v \sigma^2-\sigma_J^2=\frac{1}{p_c\big(1+\partial_p\hat \Pi_{p_c}\big)}\sigma_J^2-\sigma_J^2=\frac{\hat \Pi_{p_c}-p_c\partial_p\hat \Pi_{p_c}}{p_c\big(1+\partial_p\hat \Pi_{p_c}\big)}\sigma_J^2. 
\end{align}
By \refe{piNinfty1}, the numerator in the expression given in \refe{vsig-sigJ} can be bounded below as  
\begin{align}\lbeq{PiderivativePi}
\hat \Pi_{p_c}-p_c\partial_p\hat \Pi_{p_c}=\sum_{m=2}^{\infty}(1-m)\hat \pi_{p_c,m}&=\sum_{m=2}^{\infty}(m-1)\hat \pi_{p_c, m}^{\sss(1)}-\sum_{N=2}^{\infty}(-1)^N\sum_{m=2}^{\infty}(m-1)\hat \pi_{p_c, m}^{\sss(N)}\nn\\
&\ge \hat \pi_{p_c, 2}^{\sss(1)}-C\beta^2. 
\end{align}
Since $\hat \pi_{p_c, 2}^{\sss(1)}=p_cD^{*2}(o)$, which is $O(\beta)$ by \refe{D*2beta}, it follows that \refe{PiderivativePi} is positive for sufficiently large $L$. 

Furthermore, by \refe{ee}, we obtain
\begin{align}
\int_{\tilde S_{\tau},~\norm{l}_{\infty} \le\sqrt{t}}[\cdots]d^dl&\le \int_{\tilde S_{\tau},~\norm{l}_{\infty} \le \sqrt{t}}\Big|e^{-t\big[1-\hat J_{p_c}\big(\frac{l}{L\sqrt{t}}\big)\big]}e^{-\tau\frac{(v\sigma^2-\sigma_J^2)|l|^2}{2dL^2t}}-e^{-\frac{\sigma_J^2 |l|^2}{2dL^2}}\Big|\nn\\
&=\int_{\tilde S_{\tau},~\norm{l}_{\infty} \le\sqrt{t}}\Big|e^{-\frac{\sigma_J^2 |l|^2}{2dL^2}}e^{-\tau\frac{(v\sigma^2-\sigma_J^2)|l|^2}{2dL^2t}}e^{-t\Big[1-\hat J_{p_c}\big(\frac{l}{L\sqrt{t}}\big)-\frac{\sigma_J^2|l|^2}{2dL^2t}\Big]}-e^{-\frac{\sigma_J^2 |l|^2}{2dL^2}}\Big|\nn\\
&=\int_{\tilde S_{\tau},~\norm{l}_{\infty} \le\sqrt{t}}\Big|e^{-\frac{\sigma_J^2 |l|^2}{2dL^2}}(1-E_1)(1-E_2)-e^{-\frac{\sigma_J^2 |l|^2}{2dL^2}}\Big|\nn\\
&=\int_{\tilde S_{\tau},~\norm{l}_{\infty} \le\sqrt{t}}\Big|e^{-\frac{\sigma_J^2 |l|^2}{2dL^2}}(E_2-1)E_1-e^{-\frac{\sigma_J^2 |l|^2}{2dL^2}}E_2\Big|\nn\\
&\le\int_{\mathbb R^d}\Big\{|E_1|e^{-t\big[1-\hat J_{p_c}\big(\frac{l}{L\sqrt{t}}\big)\big]}+|E_2|e^{-\frac{\sigma_J^2 |l|^2}{2dL^2}}\Big\}\ind{\tilde S_{\tau},~\norm{l}_{\infty} \le\sqrt{t}}\lbeq{ff}
\end{align}
where 
\begin{align}
E_1=E_1(t,\tau,L)&:=1-\exp\Big(-\tau\frac{(v\sigma^2-\sigma_J^2)|l|^2}{2dL^2t}\Big),\lbeq{n}\\
E_2=E_2(t,L)&:=1-\exp\Big(-t\Big[1-\hat J_{p_c}\big(\frac{l}{L\sqrt{t}}\big)-\frac{\sigma_J^2|l|^2}{2dL^2t}\Big]\Big).\lbeq{i}
\end{align}
We note that $E_2$ is independent of $\tau$. To apply the dominated convergence theorem, we now show that the integrand in \refe{ff} converges to $0$ as $L \rightarrow \infty$ and $t \rightarrow \infty$. 

Since we have
\begin{align}
E_1\le 1-\exp\Big(-\frac{(v\sigma^2-\sigma_J^2)|l|^2}{2dL^2}\Big),
\end{align}
together with \refe{dd}, we obtain
\begin{align}\lbeq{ii}
E_1e^{-t\big[1-\hat J_{p_c}\big(\frac{l}{L\sqrt{t}}\big)\big]}\le \bigg[1-\exp\Big(-\frac{(v\sigma^2-\sigma_J^2)|l|^2}{2dL^2}\Big)\bigg]e^{-C|l|^2},
\end{align}
where \refe{vsig-sigJ} ensures that it converges to $0$ as $L \rightarrow \infty$.

Since by \refe{358}, we have
\begin{align}\lbeq{jj}
|E_2|&=\Bigg|1-\exp\Big(-\sum_{x \in \Zd}J_{p_c}(x)\frac{(l \cdot x)^2}{L^2}\int_{0}^1(1-s)(e^{\frac{isl \cdot x}{L\sqrt{t}}}-1)\Big)\Bigg|, 
\end{align}
owing to the last factor $e^{\frac{isl \cdot x}{L\sqrt{t}}}-1$ in \refe{jj}, we obtain
\begin{align}\lbeq{eeee222}
E_2e^{-\frac{\sigma_J^2 |l|^2}{2dL^2}}\overset{t \rightarrow \infty}{\rightarrow} 0.
\end{align}

Hence, by \refe{ii} and \refe{eeee222}, the integrand in \refe{ff} converges to $0$ as $L \rightarrow \infty$ and $t \rightarrow \infty$, which, together with \refe{ff} and the dominated convergence, implies that for any $\epsilon>0$, there exist $L_0$ and $T_0$ such that for all $(L,T)\ge (L_0, T_0)$, 
\begin{align}\lbeq{epep}
\Big|\int_{\tilde S_{\tau},~\norm{l}_{\infty} \le\sqrt{t}}[\cdots]d^dl\Big|\le \epsilon.
\end{align}
Therefore, for all $\epsilon$ and $(L_o,T_o)$ such that \refe{epep} holds, if we choose $(L_o,\tau_o)$ with $\tau_o\ge T_o$, then  for all $(L,\tau)\ge (L_o,\tau_o)$, we obtain
\begin{align}\lbeq{hh}
\beta\int_{\tau}^{\infty}t^{-\frac{d}{2}}\int_{\tilde S_{\tau},~\norm{l}_{\infty} \le\sqrt{t}}[\cdots]d^dldt\le \beta\epsilon \int_{\tau}^{\infty}t^{-\frac{d}{2}}dt\le C\epsilon \beta\tau^{-\frac{d-2}{2}}.
\end{align}
Thus, we obtain
\begin{align}
f(L,\tau):=\frac{\beta\int_{\tau}^{\infty}t^{-\frac{d}{2}}\int_{\tilde S_{\tau},~\norm{l}_{\infty} \le\sqrt{t}}[\cdots]d^dldt}{\beta\tau^{-\frac{d-2}{2}}}\le C\epsilon
\end{align}
which implies that $f(L,\tau)=o_{L\cap \tau}(1)$. Consequently, we conclude \refe{gg}.

\QED

\subsubsection{Proof of \refe{M1Ltau0}.}\label{lmmM1Ltau0proof}
Next, we prove \refe{M1Ltau0}. In the remainder of the paper, we frequently use the heat-kernel bound (see \cite[(1.10)]{hs05} and \cite{hs02} for details):
\begin{align}\lbeq{heatkernel}
\norm{D^{*n}}_{\infty}\le O(\beta)n^{-\frac{d}{2}}~~~(n\ge 1),&&\norm{\hat D^{n+2}}_1 \le O(\beta)(1\vee n)^{-\frac{d}{2}}~~~(n\ge 0).
\end{align}
We also use the followings for $d>4$ and for sufficiently large $L$: 
\begin{align}\lbeq{D*G}
\norm{D*G_{p_c}}_{\infty}\le O(\beta),&&\norm{D*G_{p_c}^{*2}}_{\infty}\le O(\beta),&&\norm{\hat G_{p_c}}_2^2\le C
\end{align}
which follows from $G_{p}(x)\le \delta_{o,x}+(D*G_{p})(x)$, \refe{infragauss} and \refe{heatkernel}. 

Additionally, we use \refe{G1} with $g_m=\hat \pi_{p_c,m}$, i.e., 
\begin{align}\lbeq{pimpck}
|\hat\pi_{p,m}(k)|\le C\beta m^{-\frac{d}{2}}~~~(m \ge 2).
\end{align}

First we establish the following lemma, which will be used in the proof of \refe{M1Ltau0}.
\begin{lmm}\label{dpiCG}
Let $d>4$. For sufficiently large $L$ and for $m\ge2$,
\begin{align}\lbeq{dpiCGnoab}
\Big|\int_{[-\pi,\pi]^d}\hat D(k)\hat\pi_{p_c,m}(k)\hat G_{p_c}(k)\frac{d^dk}{(2\pi)^d}\Big|\le
C\beta^2m^{-\frac{d}{2}}.
\end{align}
Additionally, for $m\ge 2$ and $n\ge1$, 
\begin{align}\lbeq{pipin1}
\int_{[-\pi,\pi]^d}\left|\hat D(k)\hat\pi_{p_c,m}(k)\hat C_{p_c,n}(k)\hat G_{p_c}(k)\right|\frac{d^dk}{(2\pi)^d}\le C\beta^{\frac{3}{2}}m^{-\frac{d}{2}}n^{-\frac{d-2}{2}}
\end{align}
\end{lmm}
\Proof{Proof of Lemma~\ref{dpiCG}.}
We begin with \refe{dpiCGnoab}. By using \refe{piNinfty1} and \refe{D*G}, we obtain
\begin{align}
\Big|\int_{[-\pi,\pi]^d}\hat D(k)\hat\pi_{p_c,m}(k)\hat G_{p_c}(k)\Big|&=\Big|\sum_{x \in \Zd}\pi_{p_c,m}(x)(D*G_{p_c})(x)\Big|\nn\\
&\le \norm{D*G_{p_c}}_{\infty}\sum_{N=1}^{\infty}\hat\pi_{p_c,m}^{\sss(N)}(0)\le C\beta^2m^{-\frac{d}{2}}.
\end{align}

For $n\ge 1$, we establish the followings. 
\begin{align}
\int_{S_n}\left|\hat D(k)\hat\pi_{p_c,m}(k)\hat C_{p_c,n}(k)\hat G_{p_c}(k)\right|&\le C\beta^2m^{-\frac{d}{2}}n^{-\frac{d-2}{2}},\lbeq{dpiCGS}\\
\int_{L_n}\left|\hat D(k)\hat\pi_{p_c,m}(k)\hat C_{p_c,n}(k)\hat G_{p_c}(k)\right|&\le C\beta^{\frac{3}{2}}m^{-\frac{d}{2}}n^{-\frac{d-2}{2}}(\log n)^{-1}.\lbeq{dpiCGL}
\end{align}
We first consider \refe{dpiCGS}. By using \refe{pimpck}, the left-hand side in \refe{dpiCGS} can be bounded as
\begin{align}\lbeq{DPCGStau}
\int_{S_n}\left|\hat D(k)\hat\pi_{p_c,m}(k)\hat C_{p_c,n}(k)\hat G_{p_c}(k)\right|&\le C\beta m^{-\frac{d}{2}}\int_{S_n}\left|\hat D(k)\hat C_{p_c,n}(k)\hat G_{p_c}(k)\right|.
\end{align}
As in the proof of \refe{M1Stau}, using \refe{localCLT} and \refe{kernelG}, we can estimate the integral in \refe{DPCGStau} as 
\begin{align}
\int_{S_n}\left|\hat D(k)\hat C_{p_c,n}(k)\hat G_{p_c}(k)\right|&\le\left[1+O\big(\beta n^{-\frac{(d-4)}{2}}\big)\right]M_{1.1}(n)+M_{1.2}(n)
\end{align}
where $M_{1.1}(n)$ and $M_{1.2}(n)$ are defined as in \refe{M11}--\refe{M12}, with $\tau$ replaced by $n$ and an absolute value applied to the integrand. Using the same approach as in the proof of \refe{M11final} and \refe{approxM12}, we can show that $M_{1.1}(n)\le C\beta n^{-\frac{d-2}{2}}$ and $M_{1.2}(n)\le C'\beta n^{-\frac{d-2+2\delta}{2}}$. Hence, we obtain \refe{dpiCGS}. 

Next, we prove \refe{dpiCGL}. By \refe{bCL1}, \refe{pimpck} and \refe{infragauss}, the left-hand side in \refe{dpiCGL} can be bounded as
\begin{align}\lbeq{dpiCGLlast}
\int_{L_n}\left|\hat D(k)\hat\pi_{p_c,m}(k)\hat C_{p_c,n}(k)\hat G_{p_c}(k)\right|&\le C\beta m^{-\frac{d}{2}}n^{-\frac{d}{2}}\int_{L_n}\big|\hat D(k)\big|\frac{1}{[1-\hat D(k)]^{3+\rho}}\nn\\
&\le C'\beta m^{-\frac{d}{2}}n^{-\frac{d-2}{2}}(\log n)^{-1}\int_{L_n}\big|\hat D(k)\big|\frac{1}{[1-\hat D(k)]^{2+\rho}}
\end{align}
where for the last inequality we use $1-\hat D(k) \ge \gamma n^{-1}\log n$. By partitioning the integral domain into $\norm{k}_{\infty}\le L^{-1}$ and $\norm{k}_{\infty}\ge L^{-1}$ and using \refe{LUL21-hatD}--\refe{Leta1-hatD}, we can bound the integral in \refe{dpiCGLlast} above by
\begin{align}\lbeq{dpiCGLint}
&c_1^{-(2+\rho)}\int_{L_n, \norm{k}_{\infty}\le L^{-1}}\frac{1}{L^{4+2\rho}|k|^{4+2\rho}}+\eta^{-(2+\rho)}\int_{L_n,\norm{k}_{\infty}\ge L^{-1}}\big|\hat D(k)\big|\nn\\
&\le C\beta \int_{\norm{l}_{\infty}\le 1}|l|^{-4-2\rho}\frac{d^dl}{(2\pi)^d}+\eta^{-(2+\rho)}\norm{\hat D^2}_1^{1/2}. 
\end{align}
Here, we apply the change of variables from $k$ to $l:=kL$ for the first term, and the Cauchy-Schwarz inequality for the second term. Since $\rho <\frac{d-4}{2}$ as mentioned below \refe{bCL2}, the integral in \refe{dpiCGLint} is bounded and from \refe{heatkernel}, $\norm{\hat D^2}_1 \le O(\beta)$. Thus, we obtain \refe{dpiCGL}. This completes the proof of Lemma~\ref{dpiCG}.\QED

\Proof{Proof of \refe{M1Ltau0}.}
By using \refe{FTrecursionGp}, the left-hand side of \refe{M1Ltau0} can be decomposed as
\begin{align}
\int_{L_{\tau}}\hat D(k)\hat C_{p_c,\tau}(k)\hat G_{p_c}(k)&=p_cM_{1,3}+M_{1,4},
\end{align}
where
\begin{align}
M_{1.3}&=\int_{L_{\tau}}\hat D(k)^2\hat C_{p_c,\tau-1}(k)\hat G_{p_c}(k),\lbeq{defM13}\\
M_{1.4}&=\sum_{m=2}^{\tau}\int_{L_{\tau}}\hat D(k)\hat\pi_{p_c,m}(k)\hat C_{p_c,\tau-m}(k)\hat G_{p_c}(k). \lbeq{M14}
\end{align}

We begin with $M_{1.3}$. To apply \refe{bCL1}--\refe{bCL2}, we decompose $M_{1.3}$ as
\begin{align}
M_{1.3}&=\int_{L_{\tau}}\hat D(k)^2\hat C_{p_c,\tau}(k)\hat G_{p_c}(k)+\int_{L_{\tau}}\hat D(k)^2\big[\hat C_{p_c,\tau-1}(k)-\hat C_{p_c,\tau}(k)\big]\hat G_{p_c}(k).
\end{align}
By applying \refe{bCL1} to the first term and \refe{bCL2} to the second term, we obtain
\begin{align}\lbeq{M131}
&C\tau^{-\frac{d}{2}}\int_{L_{\tau}}\hat D(k)^2\frac{1}{[1-\hat D(k)]^{1+\rho}}\bigg[1+\frac{1}{1-\hat D(k)}\bigg]\big|\hat G_{p_c}(k)\big|\nn\\
&\le C'\tau^{-\frac{d-2}{2}}(\log \tau)^{-1}\int_{L_{\tau}}\hat D(k)^2\frac{1}{[1-\hat D(k)]^{2+\rho}},
\end{align}
where for the last inequality we use  \refe{infragauss} and that $1-\hat D(k) \ge \gamma \tau^{-1}\log \tau$ for $k \in L_{\tau}$. As we did in \refe{dpiCGLint}, by decomposing the integral domain into $\norm{k}_{\infty} \ge L^{-1}$ and  $\norm{k}_{\infty} \le L^{-1}$, we can show that the integral in \refe{M131} can be bounded above by $C\beta$. Hence, we obtain 
\begin{align}
|M_{1.3}|\le C\beta\tau^{-\frac{d-2}{2}}(\log \tau)^{-1}.
\end{align}

Next we consider $M_{1.4}$. We begin by deriving the bound for the summand in \refe{M14}. When $m=\tau$, applying the Cauchy-Schwarz inequality and \refe{heatkernel}--\refe{pimpck}, yields
\begin{align}\lbeq{esDpiC0G}
\int_{L_{\tau}}\hat D(k)\hat\pi_{p_c,\tau}(k)\hat G_{p_c}(k)\le C\beta\tau^{-\frac{d}{2}}\int_{[-\pi,\pi]^d}\big|\hat D(k)\hat G_{p_c}(k)\big|&\le C'\beta\tau^{-\frac{d}{2}}\norm{\hat D^2}_1^{\frac{1}{2}}\norm{\hat G_{p_c}}_2\nn\\
&\le C\beta^{\frac{3}{2}}\tau^{-\frac{d}{2}}.
\end{align}
For $m<\tau$, by using $L_{\tau}\subset L_{\tau-m}$ and \refe{dpiCGL}, we obtain 
\begin{align}\lbeq{esDpiCG}
\int_{L_{\tau}}\hat D(k)\hat\pi_{p_c,m}(k)\hat C_{p_c,\tau-m}(k)\hat G_{p_c}(k)&\le \int_{L_{\tau-m}}\big|\hat D(k)\hat\pi_{p_c,m}(k)\hat C_{p_c,\tau-m}(k)\hat G_{p_c}(k)\big|\nn\\
&\le C\beta^{\frac{3}{2}}m^{-\frac{d}{2}}(\tau-m)^{-\frac{d-2}{2}}(\log (\tau-m))^{-1}.
\end{align}
Therefore, by \refe{M14} and \refe{esDpiC0G}--\refe{esDpiCG} can be bounded as 
\begin{align}\lbeq{M142}
|M_{1.4}|&\le C\beta^{\frac{3}{2}}\sum_{m=2}^{\tau}m^{-\frac{d}{2}}(\tau-m+1)^{-\frac{d-2}{2}}\nn\\
&\le C\beta^{\frac{3}{2}}\bigg\{\tau^{-\frac{d-2}{2}}\sum_{m=2}^{\lfloor(\tau+1)/2\rfloor}m^{-\frac{d}{2}}+\tau^{-\frac{d}{2}}\sum_{m=\lfloor(\tau+1)/2\rfloor +1}^{\tau}(\tau-m+1)^{-\frac{d-2}{2}}\bigg\}\nn\\
&\le C\beta^{\frac{3}{2}}\tau^{-\frac{d-2}{2}}
\end{align}
where, in the second inequality, we bound the larger of $m$ or $\tau-m+1$ below by $\tau/2$. We can then conclude the proof.
\QED

\subsection{Proof of Lemma~\ref{lmm:M2M3}}\label{derivativePiandintersectionPi1}
This section provides the proof of Lemma~\ref{lmm:M2M3}.
 \Proof{Proof of Lemma~\ref{lmm:M2M3}.}
For $M_2$, \refe{dpiCGnoab} immediately implies 
\begin{align}
|M_2|=\big|\int_{[-\pi,\pi]^d}\hat D(k)\hat \pi_{p_c,\tau}(k)\hat G_{p_c}(k)\frac{d^dk}{(2\pi)^d}\big|\le C\beta^2\tau^{-\frac{d}{2}}.
\end{align}

For $M_3$, recalling the definition of $\hat E_{p,\tau}$ from \refe{Ehatdef}, we decompose $M_{3}$ as 
\begin{align}
M_3&=M_{3.1}+M_{3.2},
\end{align}
where
\begin{align}
M_{3.1}&=\sum_{m=2}^{\tau-1}\sum_{n=\tau+1-m}^{\tau-1}\int_{[-\pi,\pi]^d}\hat D(k)\hat \pi_{p_c,m}(k)\hat C_{p_c,n}(k)\hat G_{p_c}(k),\\
M_{3.2}&=\sum_{m=\tau}^{\infty}\sum_{n=0}^{\tau-1}\int_{[-\pi,\pi]^d}\hat D(k)\hat \pi_{p_c,m}(k)\hat C_{p_c,n}(k)\hat G_{p_c}(k).\lbeq{M32}
\end{align}
By Lemma~\ref{dpiCG}, $M_{3.1}$ can be estimated similarly to \refe{M142}, yielding
\begin{align}\lbeq{M31}
|M_{3.1}|&\le C\beta^{\frac{3}{2}}\sum_{m=2}^{\tau-1}\sum_{n=\tau+1-m}^{\tau-1}m^{-\frac{d}{2}}n^{-\frac{d-2}{2}}\nn\\
&\le C\beta^{\frac{3}{2}}\bigg\{\tau^{-\frac{d-2}{2}}\sum_{m=2}^{\lfloor(\tau+1)/2\rfloor}m^{-\frac{d-2}{2}}+\tau^{-\frac{d}{2}}\sum_{m=\lfloor(\tau+1)/2\rfloor+1}^{\tau-1}\sum_{n=\tau+1-m}^{\tau-1}n^{-\frac{d-2}{2}}\bigg\}\nn\\
&\le C\beta^{\frac{3}{2}}\tau^{-\frac{d-2}{2}}.
\end{align}
We finally consider $M_{3.2}$. Since the summand in \refe{M32} can be expressed as the convolution of $\pi_{p_c,m}$ and $D*C_{p_c,n}*G_{p_c}$, by \refe{piNinfty1} and \refe{D*G}, $M_{3.2}$ can be bounded as
\begin{align}\lbeq{M322}
|M_{3.2}|\le \Big|\sum_{m=\tau}^{\infty}\sum_{y \in \Zd}\pi_{p_c,m}(y)\big(D*\sum_{n=0}^{\tau-1}C_{p_c,n}*G_{p_c}\big)(y)\Big|
&\le \norm{D*G_{p_c}^{*2}}_{\infty}\sum_{m=\tau}^{\infty}\sum_{N=1}^{\infty}\hat \pi_{p_c,m}^{\sss(N)}(0)\nn\\
&\le C\beta^2\tau^{-\frac{d-2}{2}}. 
\end{align}
Thus, combining \refe{M31}--\refe{M322}, we obtain $|M_{3}|\le C\beta^{\frac{3}{2}}\tau^{-\frac{d-2}{2}}$ as desired, completing the proof. \QED

\section{Bounds on the lace-expansion coefficients}\label{Pi2proof}
In this section, we prove Proposition~\ref{lmm:error1}. The upper bound for $\hat \pi_{p,m}^{\sss\tau}(k)$ follows from the following lemma.

\begin{lmm}\label{lmm:memorytau}
Let $d>4$. For sufficiently large $L$, the corresponding statement in \refe{piNinfty1} also holds for the memory-$\tau$ walk. Consequently, \refe{G1}--\refe{G3} hold for $m\ge 1$ in the memory-$\tau$ walk as well.  
\end{lmm}
\Proof{Proof of Lemma~\ref{lmm:memorytau}.}
As noted in Section~\ref{results}, van der Hofstad and Slade established in \cite[Proposition~4.1]{hGs03} that Assumption G holds for the self-avoiding walk by proving \refe{piNinfty1}, given \refe{127hs03} as a hypothesis. To establish \refe{piNinfty1} using \refe{127hs03}, they proceed by induction on $N$ for $\pi_{p,m}^{\sss(N)}$. In the proof, they decompose the path considered for $\pi_{p,m}^{\sss(N)}$ into self-avoiding subpaths according to the time structure given by the graph of a lace in $\mathcal L^{\sss(N)}$ and then apply the assumptions from \refe{127hs03} to each subpath to derive \refe{piNinfty1}. Consequently, \refe{127hs03} are shown to hold for all $m\ge1$ for the self-avoiding walk by induction.

For the memory-$\tau$ walk, we can apply the same method as for the self-avoiding walk, decomposing the path considered for $\pi_{p,m}^{\sss\tau, (N)}$ into memory-$\tau$ self-avoiding subpaths based on the time structure determined by a lace in $\mathcal L_{\tau}^{\sss(N)}$. It is important to note that in the case of the memory-$\tau$ walk, the subpaths obtained from the decomposition of $\pi_{p,m}^{\sss\tau, (N)}$ are all shorter than $\tau$ in length. Consequently, these subpaths are identical to self-avoiding paths, enabling us to apply \refe{127hs03} for the self-avoiding walk as a fact, rather than a hypothesis, to each subpath, even in the case of the memory-$\tau$ walk. For this reason, we can establish \refe{piNinfty1} for the memory-$\tau$ walk using the same approach. For the proof of \refe{piNinfty1}, see \cite[Section~5]{hGs03}. \QED

\Proof{Proof of Proposition~\ref{lmm:error1}.}First, we observe that by \refe{pinpitau}, it follows that for $N \ge2$
\begin{align}\lbeq{aaa}
\hat \Pi_{p_c^{\sss\tau}}^{\sss(N)}-\hat \Pi_{p_c^{\sss\tau}}^{\sss\tau,(N)}=\sum_{n=\tau+1}^{\infty}\big[\hat\pi_{p_c^{\sss\tau}, n}^{\sss(N)}-\hat\pi_{p_c^{\tau}, n}^{\sss\tau, (N)}\big].
\end{align}
Thus, by \refe{piNinfty1} and Lemma~\ref{lmm:memorytau}, we obtain 
\begin{align}\lbeq{finalPi2Pi2}
C'\beta^2 \tau^{-\frac{d-2}{2}}\ge \sum_{N=2}^{\infty}{(-1)}^N\left[\hat \Pi_{p_c^{\sss\tau}}^{\sss(N)}
-\hat \Pi_{p_c^{\tau}}^{\sss\tau,(N)}\right]\ge \sum_{n=\tau+1}^{\infty}\big[\hat\pi_{p_c^{\sss\tau}, n}^{\sss(2)}-\hat\pi_{p_c^{\tau}, n}^{\sss\tau, (2)}\big]-C\beta^3\tau^{-\frac{d-2}{2}}.
\end{align}
To complete the proof, it is therefore enough to show the existence of a positive constant $K$ such that 
\begin{align}
\sum_{n=\tau+1}^{\infty}\big[\hat\pi_{p_c^{\sss\tau}, n}^{\sss(2)}-\hat\pi_{p_c^{\tau}, n}^{\sss\tau, (2)}\big]\ge K\beta^2\tau^{-\frac{d-2}{2}}. 
\end{align}
We note that since for the memory-$\tau$ walk, the length of an edge in a lace is at most $\tau$, $\pi_{p, n}^{\sss\tau, (2)}(x)=0$ for $n \ge 2\tau$. Therefore, recalling the definitions given in \refe{mathcalJN}--\refe{pipnalpha}, we obtain
\begin{align}\lbeq{pi2pi2tau}
\sum_{n=\tau+1}^{\infty}\big[\hat\pi_{p_c^{\sss\tau}, n}^{\sss(2)}-\hat\pi_{p_c^{\tau}, n}^{\sss\tau, (2)}\big]&=\sum_{n=2\tau}^{\infty}\hat\pi_{p_c^{\sss\tau}, n}^{\sss(2)}+\sum_{n=\tau+1}^{2\tau-1}\sum_{x \in \Zd}\sum_{\omega \in \mathcal W_n(o,x)}W_{p_c^{\tau}}(\omega)\big[\mathcal J^{\sss(2)}[0,n]-\mathcal J_{\tau}^{\sss(2)}[0,n]\big]\nn\\
&\ge \sum_{n=2\tau}^{\infty}\hat\pi_{p_c^{\sss\tau}, n}^{\sss(2)}-\sum_{n=\tau+1}^{2\tau-1}\sum_{x \in \Zd}\sum_{\omega \in \mathcal W_n(o,x)}W_{p_c^{\sss\tau}}(\omega)\tilde{\mathcal J}_{\tau}^{\sss(2)}[0,n]
\end{align}
where 
\begin{align}
\tilde{\mathcal J}_{\tau}^{\sss(2)}[0,n]=\sum_{L\in \mathcal L_{\tau}^{\sss(2)}[a,b]}\prod_{st \in L}(-\mathcal U_{st})\prod_{s't' \in \mathcal C_{\tau}(L)}(1+\mathcal U_{s't'})\Big[1-\prod_{s't' \in \mathcal C_{\tau}(L)\backslash \mathcal C(L)}(1+\mathcal U_{s't'})\Big].
\end{align}
For the second line of \refe{pi2pi2tau}, we use the inclusions $\mathcal L_{\tau}^{\sss(2)}[a,b]\subset \mathcal L^{\sss(2)}[a,b]$ and $ \mathcal C_{\tau}(L)\subset \mathcal C(L)$, and neglect the contribution from the summation over $L \in \mathcal L^{\sss(2)}[a,b]\setminus\mathcal L_{\tau}^{\sss(2)}[a,b]$, because it is non-negative due to $L$ containing exactly two edges.

\begin{figure}[t]
\begin{center}
\includegraphics[scale=0.6]{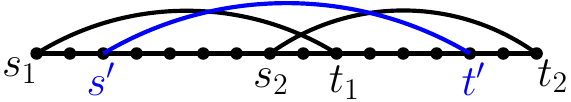}\hskip15mm\includegraphics[scale=0.6]{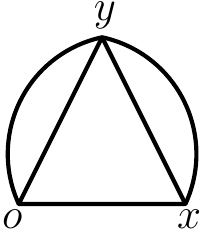}
\end{center}
\caption{A sample lace $L$ in $\mathcal L^{\sss(2)}_{\tau}$ and a compatible edge $s't'$ in $\mathcal C_{\tau}(L)\backslash \mathcal C(L)$. The black edges are shorter than $\tau$, while the blue edge is longer than $\tau+1$ in length. The right figure provides a diagrammatic representation of the left graph, where each black line denotes a self-avoiding subpath of $\omega$. In particular, $\omega_{s_1}=\omega_{t_1}=o$, $\omega_{s_2}=\omega_{t_2}=x$ and $\omega_{s'}=\omega_{t'}=y$. These subpaths exhibit mutual avoidance effects. }
\label{fig:Lace2}
\end{figure}

First we derive an upper bound for the second term of \refe{pi2pi2tau}. Observe that $\mathcal C_{\tau}(L)\backslash \mathcal C(L)$ consists of compatible edges with $L$ whose lengths are at least $\tau+1$. Since the length of edge is at most $\tau$, the endpoints of any edge $\in \mathcal C_{\tau}(L)\backslash \mathcal C(L)$ cannot lie within the interval of the same edge in a lace (see Figure~\ref{fig:Lace2}). 
Thus, since we have
\begin{align}
1-\prod_{s't' \in \mathcal C_{\tau}(L)\backslash \mathcal C(L)}(1+\mathcal U_{s't'})\le \sum_{s't'\in \mathcal C_{\tau}(L)\backslash \mathcal C(L)}\ind{\omega_{s'}=\omega_{t'}},
\end{align}
and by neglecting the mutual avoidance among subpaths, the second term in \refe{pi2pi2tau} can be bounded as
\begin{align}\lbeq{CCCCC}
&\sum_{n=\tau+1}^{2\tau-1}\sum_{x \in \Zd}\sum_{\omega \in \mathcal W_n(o,x)}W_{p_c^{\sss\tau}}(\omega)\tilde{\mathcal J}_{\tau}^{\sss(2)}[0,n]\nn\\
&\le \sum_{n=\tau+1}^{2\tau-1}\sum_{x,y \in \Zd}\sum_{\substack{n_1,\cdots, n_5\\ (\star) }}C_{n_1}(y)C_{n_2}(y)C_{n_3}(x)C_{n_4}(x-y)C_{n_5}(x-y). 
\end{align}
In the above, we omit $p_c^{\sss\tau}$ in the indices of $C_{n_i}$, and the condition $(\star)$ signifies that $ n_1+\cdots+n_5=n$ and $1\le n_i\le \tau$ for $1\le i\le 5$. Note that each $C_{n_i}$ is defined for the self-avoiding walk, as the length of each subpath is less than $\tau$, making each subpath a self-avoiding walk rather than a memory-$\tau$ walk. By bounding $C_{n_i}$ with the largest index using $L_{\infty}$-norm and applying \refe{127hs03} with $f_m=\hat C_m$, \refe{CCCCC} can be bounded above by 
\begin{align}
C\beta \sum_{n=\tau+1}^{2\tau-1}n^{-\frac{d}{2}}\norm{D*G_{p_c^{\sss\tau}}^{*2}}_{\infty}^2\le C'\beta^3\tau^{-\frac{d-2}{2}}.
\end{align}
Here, we use the bound $C_{n}(x) \le (D*C_{n-1})(x)$ for any $n\ge 1$, which results from disregarding mutual avoidance, together with \refe{D*G} to obtain the inequality. 

Next, we consider the first term of \refe{pi2pi2tau}. By \refe{pctau-pcinf}, Proposition~\ref{lmm:error0} and the upper bound from \refe{finalPi2Pi2}, we have $p_c-p_c^{\sss\tau}\le C\beta\tau^{-\frac{d-2}{2}}$. Thus, by applying \refe{piNinfty1}, we obtain
\begin{align}
\sum_{n=2\tau}^{\infty}\hat\pi_{p_c^{\sss\tau}, n}^{\sss(2)}\ge \sum_{n=2\tau}^{\infty}\hat\pi_{p_c, n}^{\sss(2)}-(p_c-p_c^{\sss\tau})\sum_{n=2\tau}^{\infty}\partial_p\hat\pi_{p_c, n}^{\sss(2)}\ge \sum_{n=2\tau}^{\infty}\hat\pi_{p_c, n}^{\sss(2)}-C\beta^3\tau^{-(d-3)}
\end{align}
To obtain a lower bound of $\hat\pi_{p_c, n}^{\sss(2)}$, we focus on a specific lace: $L=\{0t_1, s_2n\}$, where $t_1=n-1$ and $s_2=n-2$. Then, we obtain
\begin{align}
\hat\pi_{p_c, n}^{\sss(2)}&\ge \sum_{x \in \Zd}\sum_{\omega \in \mathcal W_n(o,x)}W_{p_c}(\omega)\mathcal U_{0(n-1)}\mathcal U_{(n-2)n}\prod_{s't' \in \mathcal C(L)}(1+\mathcal U_{s't'})\nn\\
&=p_c^2\sum_{x\in \Zd}D(x)^2C_{p_c,n-2}(x). 
\end{align}
Since Lemma~\ref{lmm:rewriteCbiggertau} provides
\begin{align}
\sum_{n=2\tau}^{\infty}C_{p_c,n-2}(x)=(G_{p_c}*C_{p_c,2\tau-2})(x)-(G*\pi_{p_c,2\tau-2})(x)+(G_{p_c}*E_{p_c,2\tau-2})(x),
\end{align}
we obtain
\begin{align}
\sum_{n=2\tau}^{\infty}\hat\pi_{p_c, n}^{\sss(2)}&\ge p_c^2\sum_{x\in \Zd}D(x)^2(G_{p_c}*C_{p_c,2\tau-2})(x)-p_c^2C_2\beta\big\{|M_2(2\tau-2)|+|M_3(2\tau-2)|\big\}\nn\\
&\ge p_c^2\sum_{x\in \Zd}D(x)^2(G_{p_c}*C_{p_c,2\tau-2})(x)-C\beta^{\frac{5}{2}}\tau^{-\frac{d-2}{2}}\big[\beta^{\frac{1}{2}}\tau^{-1}+1\big]
\end{align}
where $M_2(2\tau-2)$ and $M_3(2\tau-2)$ are defined in \refe{defM2}--\refe{defM3}, respectively, with $\tau$ replaced by $2\tau-2$. Here, the first inequality follows from \refe{sigmabound}, while the second inequality follows from Lemma~\ref{lmm:M2M3}.

In the rest of the proof, we aim to establish that there exisits a positive constant $C$ such that 
\begin{align}\lbeq{D2GCn-2}
\sum_{x\in \Zd}D(x)^2(G_{p_c}*C_{p_c,2\tau-2})(x)\ge C\beta^2\tau^{-\frac{d-2}{2}}. 
\end{align}
The left-hand side of \refe{D2GCn-2} can be rewritten in terms of the Fourier transform as
\begin{align}\lbeq{D2StauandLtau}
\int_{S_{2\tau-2}}\hat {D^2}(k)\hat G_{p_c}(k)\hat C_{p_c,2\tau-2}(k)+\int_{L_{2\tau-2}}\hat {D^2}(k)\hat G_{p_c}(k)\hat C_{p_c,2\tau-2}(k)
\end{align}
noting that $\hat {D^2}(k)=\sum_{x}D(x)^2e^{ik\cdot x}$. 

We first show that the the second term acts as an error term with respect to either $\beta$ or $\tau$, following a proof similar to that of \refe{M1Ltau0}. Using \refe{FTrecursionGp}, the second term of \refe{D2StauandLtau} can be decomposed into $\tilde M_{1.3}(2\tau-1)$ and $\tilde M_{1.4}(2\tau-1)$, which  are defined in \refe{defM13} and \refe{M14} respectively, with $\tau$ replaced by $2\tau-1$ and a single occurence of $\hat D(k)$ replaced by $\hat {D^2}(k)$. For $\tilde M_{1.3}(2\tau-1)$, as in \refe{M131}, we obtain
\begin{align}\lbeq{tildeM13}
|\tilde M_{1.3}(2\tau-1)|\le C'\tau^{-\frac{d-2}{2}}(\log \tau)^{-1}\int_{L_{2\tau-2}}|\hat {D^2}(k)\hat D(k)|\frac{1}{[1-\hat D(k)]^{2+\rho}}.
\end{align}
Noting that $|\hat {D^2}(k)|\le O(\beta)$, $\norm{(\hat {D^2})^2}_1\le O(\beta^{3})$, we apply the Cauchy-Schwarz inequality to bound the integral in \refe{tildeM13} by $C\beta^2$. Similarly, for $\tilde M_{1.4}(2\tau-1)$, we obatin
\begin{align}
\Big|\int_{L_{\tau}}\hat D(k)\hat\pi_{p_c,\tau}(k)\hat G_{p_c}(k)\Big|&\le C\beta^{\frac{5}{2}}\tau^{-\frac{d}{2}}\\
\Big|\int_{L_{2\tau-2}}\hat {D^2}(k)\hat\pi_{p_c,m}(k)\hat C_{p_c,2\tau-2-m}(k)\hat G_{p_c}(k)\Big|&\le C\beta^{\frac{5}{2}}m^{-\frac{d}{2}}\tau^{-\frac{d-2}{2}}(\log \tau)^{-1}.
\end{align}
Therefore, as in \refe{M142}, $\tilde M_{1.4}(2\tau-1)$ can be bounded above by $C\beta^{\frac{5}{2}}\tau^{-\frac{d-2}{2}}$, and together with the upper bound of $\tilde M_{1.3}(2\tau-1)$, this implies that the second term of \refe{D2StauandLtau} can be bounded as
\begin{align}\lbeq{LtauN2}
\int_{L_{2\tau-2}}\hat {D^2}(k)\hat G_{p_c}(k)\hat C_{p_c,2\tau-2}(k)\le C\beta^2\tau^{-\frac{d-2}{2}}(\log \tau)^{-1}+C'\beta^{\frac{5}{2}}\tau^{-\frac{d-2}{2}}.
\end{align}

Finally, we consider the first term of \refe{D2StauandLtau}, which can be decomposed as 
\begin{align}\lbeq{hdhshv}
D^{*2}(o)\int_{S_{2\tau-2}}\hat G_{p_c}(k)\hat C_{p_c,2\tau-2}(k)+\int_{S_{2\tau-2}}\big[\hat {D^2}(k)-D^{*2}(o)\big]\hat G_{p_c}(k)\hat C_{p_c,2\tau-2}(k). 
\end{align}
For the first term, we apply the same argument as in \refe{M1Stau}: the integral in the first term of \refe{hdhshv} corresponds to \refe{M1Stau} with $\tau$ replaced by $2\tau-2$ and with $\hat D(k)$ omitted from the integrand. Accordingly, we split it into two parts, analogous $M_{1.1}$ and $M_{1.2}$ in the proof of \refe{M1Stau}, but without $\hat D(k)$. Then, together with \refe{D*2beta}, we conclude that there exists a positive constant C such that 
\begin{align}\lbeq{mainN2}
D^{*2}(o)\int_{S_{2\tau-2}}\hat G_{p_c}(k)\hat C_{p_c,2\tau-2}(k)\ge C\beta^2\tau^{-\frac{d-2}{2}}.
\end{align}

Thus, the remaining task is to show
\begin{align}\lbeq{hatDDo}
\Big|\int_{S_{2\tau-2}}\big[\hat {D^2}(k)-D^{*2}(o)\big]\hat G_{p_c}(k)\hat C_{p_c,2\tau-2}(k)\Big|\le C\beta^2\tau^{-\frac{d}{2}}\big[1+\tau^{-\delta}\big].
\end{align}
Since by using $D^{*2}(o)=\sum_{x}D(x)^2$ and \refe{taylor}, we have
\begin{align}
|\hat {D^2}(k)-D^{*2}(o)|\le |k|^2\sum_{x \in \Lambda}D(x)^2|x|^2\le C_2\beta\sigma^2|k|^2,
\end{align}
the left-hand side in \refe{hatDDo} can be bounded above by 
\begin{align}
 C_2\beta\sigma^2\int_{S_{2\tau-2}}|k|^2\big|\hat G_{p_c}(k)\hat C_{p_c,2\tau-2}(k)\big|&\le C\beta\sigma^2\int_{S_{2\tau-2}}\frac{|k|^2}{1-\hat D(k)}\big|\hat C_{p_c,2\tau-2}(k)\big|\nn\\
&\le C\beta\sigma^2\left\{\left[1+O\big(\beta \tau^{-\frac{(d-4)}{2}}\big)\right]M_{4.1}+M_{4.2}\right\}
\end{align}
where 
\begin{align}
M_{4.1}=\int_{S_{2\tau-2}}\frac{|k|^2}{1-\hat D(k)}e^{-(\tau-1)\frac{v\sigma^2|k|^2}{d}},&&M_{4.2}=\int_{S_{2\tau-2}}\frac{O(L^2|k|^4\tau^{1-\delta})}{1-\hat D(k)}e^{-(\tau-1)\frac{v\sigma^2|k|^2}{d}}.
\end{align}
In the above, we use \refe{localCLT} and \refe{infragauss}. For both $M_{4.1}$ and $M_{4.2}$, we use a similar estimate as in \refe{abM12} to obtain
\begin{align}
M_{4.1}\le C\beta L^{-2}\tau^{-\frac{d}{2}},&&M_{4.2}\le C\beta L^{-2}\tau^{-\frac{d}{2}-\delta},
\end{align}
which leads to \refe{hatDDo}, since $\sigma^2L^{-2}\le C_1$ by \refe{sigmabound}. Thus, by combining \refe{LtauN2}, \refe{mainN2}--\refe{hatDDo}, we obtain \refe{D2GCn-2}, which completes the proof. 
 \QED

\section{Notes}\label{notes}
In these notes, we provide remarks on the application of Theorem~\ref{thm:main} to other models.  The proof of Theorem~\ref{thm:main} heavily relies on Theorem~\ref{thm:local}, which was established by van der Hofstad and Slade. In other words, Theorem~\ref{thm:main} appears to extend to models for which a local central limit theorem of the type results given in Theorem~\ref{thm:local} is known. 
\begin{itemize}
\item In \cite{hhs98}, van der Hofstad, den Hollander and Slade establish the counterpart of \refe{localCLT} for the nearest-neighbor weakly self-avoiding walk, where $D(x)=\frac1{2d}\ind{\norm{x}_1=1}$ and $K[a,b]$ defined in \refe{defust} is replaced by $\prod_{a \le s<t\le b }(1+\lambda_{st}\mathcal U_{st})$ with $\lambda_{st}=\lambda_{st}(\beta,p)=1-e^{-\frac{\beta}{|s-t|^p}}$ for $p\in \mathbb R$ and $\beta\ge0$. Here, $\beta$ is used as a distinct parameter from the one defined in \refe{defbeta}. Their results hold for sufficiently small $\beta$, subject to the following conditions on $p$: for $d>4,~p\ge0$, and for $d\le4,~p>(4-d)/2$. We can extend our findings to the weakly-self-avoiding walk by applying their result from \cite{hhs98} in place of \refe{localCLT} and following the same approach as in Theorem~\ref{thm:main}. 

\item In \cite{H11}, Heydenreich establishes scaling limit-type results for the long-range self-avoiding walk above its upper critical dimension. For the long-range self-avoiding walk, $D$ is defined as \refe{defsigmah} with the conditions replaced by the following: there exists a positive constant $c_h$ such that $c_hh(x)|x|^{d+\alpha}\rightarrow 1$ as $|x|\rightarrow \infty$, with $\alpha>0$. In \cite{H11}, Heydenreich shows that with suitable scaling, the long-range self-avoiding walk in $d>2(\alpha \wedge 2)$ converges to Brownian motion for $\alpha \ge 2$, and to $\alpha$-stable L\'evy motion for $\alpha<2$. In Theorem~1.2 in \cite{H11}, the $n\rightarrow \infty$ limit is taken, but the expression given in (2.17) seems useful as an analogue of Theorem~\ref{thm:local} with some modification. Since we have the error term in (2.17) that is $O(n^{-\epsilon})$ for $\epsilon \in (0, (\frac{d}{\alpha\wedge 2}-2)\wedge 1)$, we cannot yet assert that Theorem~\ref{thm:main} extends to the long-range model. However, the author believes an extension may be achievable by adapting the results in \cite{H11}.

\item For the nearest-neighbor model in sufficiently high dimensions, Madras and Slade establish in \cite[Lemma~6.1.3]{ms93} that there exists a positive constant $B$ such that $\sup_{x}C_{p_c,n}(x)\le Bn^{-d/2}$. This represents the best estimate currently available among the local central limit theorem-type results for the nearest-neighbor self-avoiding walk. 
Using their results, it seems possible that the convergence rate of $p_c^{\sss\tau}$ to $p_c$ for the nearest-neighbor is of order $\tau^{-(d-2)/2}$. However, determining the proportionality constant of the leading term appears challenging, because the estimate given in \cite[Lemma~6.1.3]{ms93} is only an upper bound. In conclusion, to extend Theorem~\ref{thm:main} to the nearest-neighbor model, we must first establish an analogue of Theorem~\ref{thm:local}  for the nearest-neighbor case. 
\end{itemize}

\section*{Acknowledgements}The author thanks the anonymous referees for carefully reading the manuscript and providing helpful comments. The author is also grateful to Fumihiko Nakano, Akira Sakai and Gordon Slade for comments on an earlier version of this work.



\begin{thebibliography}{99}
\bibitem{bhk18}E. Bolthausen, R. van der Hofstad and G. Kozma.
\newblock Lace expansion for dummies.
\newblock \emph{Ann. Inst. Henri Poincar\'e. Probab. Stat.}, \textbf{54} (2018): 141--153.

\bibitem{bs85}D.C. Brydges and T. Spencer.
\newblock Self-avoiding walk in 5 or more dimensions.
\newblock \emph{Commun. Math. Phys.}, \textbf{97} (1985): 125--148.



\bibitem{cfs59}M.E. Fisher and M.F. Sykes.
\newblock  Excluded-volume problem and the Ising model of ferromagnetism.
\newblock \emph{Phys. Rev.}, \textbf{43} (1959): 45--58.

\bibitem{g10}B.T. Graham.
\newblock  Borel type bounds for the self-avoiding walk connective constant.
\newblock \emph{J. Phys. A: Math. Theor.}, \textbf{43} (2010): no.23,~235001,~13pp.

\bibitem{h08}T. Hara.
\newblock Decay of correlations in nearest-neighbor self-avoiding walk, percolation, lattice trees and animals. 
\newblock \emph{Ann. Probab.}, \textbf{36} (2008): 530--593.


\bibitem{hhs03}T. Hara, R. van der Hofstad and G. Slade.
\newblock Critical two-point functions and the lace expansion for spread-out 
high-dimensional percolation and related models.
\newblock \emph{Ann. Probab.}, \textbf{31} (2003): 349--408.

\bibitem{hs90a}T. Hara and G. Slade.
\newblock Mean-field critical behaviour for percolation in high dimensions.
\newblock \emph{Commun. Math. Phys.}, \textbf{128} (1990): 333--391.

\bibitem{hs90b}T. Hara and G. Slade.
\newblock On the upper critical dimension of lattice trees and lattice animals.
\newblock \emph{J. Stat. Phys.}, \textbf{59} (1990): 1469--1510.

\bibitem{hs92a}T. Hara and G. Slade.
\newblock Self-avoiding walk in five or more dimensions.~I. The critical behaviour.
\newblock \emph{Commun. Math. Phys.}, \textbf{147} (1992): 101--136.


\bibitem{hs92b}T. Hara and G. Slade.
\newblock The lace expansion for self-avoiding walk in five or more dimensions. 
\newblock \emph{Rev. Math. Phys.}, \textbf{4} (1992): 235--327.


\bibitem{hs92}T. Hara and G. Slade.
\newblock The number and size of branched polymers in high dimensions. 
\newblock \emph{J. Stat. Phys.}, \textbf{67} (1992): 1009--1038.

\bibitem{hs93}T. Hara and G. Slade.
\newblock The self-avoiding walk and percolation critical points in high dimensions. 
\newblock \emph{Comb. Probab. Comput.}, \textbf{4} (1995): 197--215.


\bibitem{hss93}T. Hara, G. Slade and Alan D.Sokal.
\newblock New lower bounds on the self-avoiding-walk connective constant. 
\newblock \emph{J. Stat. Phys.}, \textbf{72} (1993): 479--517.


\bibitem{H11}M. Heydenreich
\newblock Long-range self-avoiding walk converges to $\alpha$-stable processes. 
\newblock \emph{Ann. Inst. Henri Poincar\'e. Probab. Stat.}, \textbf{47} (2011): 20--42.


\bibitem{hhs98}R. van der Hofstad, F. den Hollander and G. Slade
\newblock A new inductive approach to the lace expansion for self-avoiding walks.
\newblock \emph{Probab. Theory Relat. Fields}, \textbf{111} (1998): 253--286.


\bibitem{hs04}R. van der Hofstad and A. Sakai.
\newblock Gaussian scaling for the critical spread-out contact process above 
the upper critical dimension.
\newblock \emph{Electron. J. Probab.}, \textbf{9} (2004): 710--769.

\bibitem{hs05}R. van der Hofstad and A. Sakai.
\newblock Critical points for spread-out self-avoiding walk, percolation and 
the contact process above the upper critical dimensions.  
\newblock \emph{Probab. Theory Relat. Fields}, \textbf{132} (2005): 438--470.

\bibitem{hGs03}R. van der Hofstad and G. Slade.
\newblock The lace expansion on a tree with application to networks of self-avoiding walks.  
\newblock \emph{Adv. in Appl. Math.}, \textbf{30} (2003): 471--528.


\bibitem{hs02}R. van der Hofstad and G. Slade.
\newblock A generalised inductive approach to the lace expansion.
\newblock \emph{Probab. Theory Relat. Fields}, \textbf{122} (2002): 389--430.


\bibitem{k67}D.A. Klarner.
\newblock Cell growth problems.
\newblock \emph{Canad. J. Math.}, \textbf{19} (1967): 851--863. 
 
\bibitem{k81}D.J. Klein.
\newblock Rigorous results for branched polymer models with excluded volume.
\newblock {\it J. Chem. Phys.}, {\bf 75} (1981): 5186--5189.

\bibitem{k64}H. Kesten.
\newblock On the number of self-avoiding walks.~II
\newblock {\it J. Math. Phys.}, {\bf 5} (1964): no.8, 1128--1137.


\bibitem{ks22}N. Kawamoto and A. Sakai.
\newblock Spread-out limit of the critical points for lattice trees and lattice animals in dimensions $d>8$.
\newblock {\it Comb. Probab. Comput.}, {\bf 33} (2024): 238--269.


\bibitem{ms93}N. Madras and G. Slade.
\newblock \emph{The Self-Avoiding Walk}, Birkh\"auser, Boston, MA,  Boston, Inc., Boston, (1993). xiv+425pp.

\bibitem{ms11}Y. M. Miranda and G. Slade.
\newblock The growth constants of lattice trees and lattice animals in high dimensions.
\newblock {\it Elect. Comm. in Probab.}, {\bf 16} (2011): 129--136.


\bibitem{ms13}Y. M. Miranda and G. Slade.
\newblock Expansion in high dimension for the growth constants of lattice trees and lattice animals.
\newblock {\it Comb. Probab. Comput.}, {\bf 22} (2013): 527--565.


\bibitem{n98}J. Noonan.
\newblock New upper bounds for the connective constants of self-avoiding walks.
\newblock {\it J. Stat. Phys.}, {\bf 91} (1998): 871--888.


\bibitem{p94}M.D. Penrose.
\newblock Self-avoiding walks and trees in spread-out lattices.
\newblock {\it J. Stat. Phys.}, {\bf 77} (1994): 3--15.

\bibitem{u02}D. Ueltschi.
\newblock Self-avoiding walk with attractive interactions.
\newblock \emph{Probab. Theory Relat. Fields}, {\bf 124} (2002): 189--203.


\bibitem{s87}G. Slade.
\newblock The diffusion of self-avoiding random walk in high dimensions.
\newblock \emph{Commun. Math. Phys.}, \textbf{110} (1987): 661--683.

\bibitem{g004}G. Slade.
\newblock \emph{The Lace Expansion and its Applications}, Springer, (2006). Lecture Notes in
 Mathematics Vol. 1879. Ecole d’Eté de Probabilités de Saint–Flour XXXIV–2004. 


\end{thebibliography}
\end{document}